\documentclass{article}

\usepackage{titling}
\usepackage{graphicx}
\usepackage{amsfonts}
\usepackage{amsmath}
\usepackage{amsthm}
\usepackage{algorithm}
\usepackage{color}
\usepackage{hyperref}
\usepackage{cite}
\usepackage{refcheck}

\newcommand{\I}{\mathrm{i}\,}
\newtheorem{theorem}{Theorem}

\theoremstyle{remark}

\title{Confluent Vandermonde with Arnoldi}

\author{Qiang Niu \qquad Hui Zhang$^*$ \qquad Youzhou Zhou\\ \small{School of Mathematics and
    Physics, Xi'an Jiaotong-Liverpool University, Ren'ai Rd. 111, Suzhou 215123, China}\\
  \small{$^*$ Corresponding author. Emails: \{qiang.niu, hui.zhang, youzhou.zhou@xjtlu.edu.cn\}}}

\date{\small{\today}}

\begin{document}
\maketitle

\begin{abstract}
  In this note, we extend the {Vandermonde with Arnoldi method} recently advocated by P. D. Brubeck,
  Y. Nakatsukasa and L. N. Trefethen [SIAM Review, 63 (2021) 405-415] to dealing with the confluent
  Vandermonde matrix. To apply the Arnoldi process, it is critical to find a Krylov subspace which
  generates the column space of the confluent Vandermonde matrix. A theorem is established for such
  Krylov subspaces for any order derivatives. This enables us to compute the derivatives of high
  degree polynomials to high precision. It also makes many applications involving derivatives
  possible, as illustrated by numerical examples. We note that one of the approaches orthogonalizes
  only the function values and is equivalent to the formula given by P. D. Brubeck and L. N.
  Trefethen [SIAM J. Sci. Comput., 44 (2022) A1205–A1226]. The other approach orthogonalizes the
  Hermite data. About which approach is preferable to another, we made the comparison, and the
  result is problem dependent.
\end{abstract}



\section{Introduction.}

The recent work \cite{va2021} presents a linear algebra method for the polynomial interpolation and
fitting on arbitrary grids. It does not need the prior knowledge of a well-conditioned polynomial
basis, and works easily on two intervals and subsets of complex plane, for example. The method
starts with the ill-conditioned Vandermonde system
\[
  \sum_{k=0}^nx_j^kc_k = f(x_j),\quad j=1,\ldots,m
\]
such that $f(x)\approx \sum_{k=0}^n x^kc_k$, and transforms it into a much more well-conditioned
system
\[
  \sum_{k=0}^nq_{jk}d_k = f(x_j),\quad j=1,\ldots,m 
\]
such that $f(s_j)\approx \sum_{k=0}^nw_{jk}d_k$. In matrix form, the Vandermonde system
$V\mathbf{c}=\mathbf{f}$, i.e.,
\[
  \begin{pmatrix}
    1 & x_1 & x_1^2 & \cdots & x_1^n\\
    1 & x_2 & x_2^2 & \cdots & x_2^n\\
    \vdots & \vdots & \vdots &\cdots & \vdots\\
    1 & x_m & x_m^2 & \cdots & x_m^n
  \end{pmatrix}
  \begin{pmatrix}
    c_0\\ c_1\\ \vdots \\ c_n
  \end{pmatrix}=
  \begin{pmatrix}
    f(x_1)\\ f(x_2)\\ \vdots\\ f(x_m)
  \end{pmatrix}
\]
is transformed into $Q\mathbf{d}=\mathbf{f}$ where $Q=(q_{jk})$ is a matrix and $\mathbf{d}=(d_k)$
is a column vector, and the approximation at $s_j$, $j=1,\ldots,M$ is evaluated as
$(f(s_j))\approx W\mathbf{d}$ where $W=(w_{jk})$. The columns of $Q$ form an orthogonal basis of the
column space of $V$. To construct $Q$, it is noted that $\mathrm{Col}\,V$ is just the Krylov
subspace $\mathrm{span}\{\mathbf{e}, X\mathbf{e}, X^2\mathbf{e}, \ldots, X^n\mathbf{e}\}$ with
$X=\mathrm{diag}(x_1,x_2,\ldots,x_m)$, $\mathbf{e}=(1,\ldots,1)^T$. So, the Arnoldi process for the
Krylov subspace is exploited to generate the columns of $Q$ as well as a Hessenberg matrix $H$ such
that $XQ_{-} = QH$ (here $Q_{-}$ is obtained by deleting the last column of $Q$). The evaluation
matrix $W$ has the first column all one's and is constructed column by column such that $SW_{-}=WH$,
where $S=\mathrm{diag}(s_1,\ldots,s_M)$ and $W_{-}$ is $W$ without the $(n+1)$-th column.

Now, we would like to employ the same idea for the Hermite interpolation and fitting. We start with
the confluent Vandermonde system
\[
  \sum_{k=0}^nx_j^kc_k = f(x_j),\quad \sum_{k=0}^nkx_j^{k-1}c_k = f'(x_j), \quad j=1,\ldots,m
\]
or in matrix form (with $V$ as before)
\[
  \begin{pmatrix}
    V\\ V_{'}
  \end{pmatrix}\mathbf{c}=
  \begin{pmatrix}
    \mathbf{f_{\phantom{'}}}\\ \mathbf{f}_{'}
  \end{pmatrix},\quad V_{'}=
  \begin{pmatrix}
    0 & 1 & 2x_1 & 3x_1^2 & \ldots & nx_1^{n-1}\\
    0 & 1 & 2x_2 & 3x_2^2 & \ldots & nx_2^{n-1}\\
    \vdots & \vdots & \vdots &\cdots & \vdots\\
    0 & 1 & 2x_m & 3x_m^2 & \ldots & nx_m^{n-1}
  \end{pmatrix},\quad\mathbf{f}_{'}=
  \begin{pmatrix}
    f'(x_1)\\ f'(x_2)\\ \vdots\\ f'(x_m)
  \end{pmatrix}.
\]
We notice that the column space of $[V; V_{'}]$ (semicolon means as in MATLAB) is a Krylov subspace:
\[
  \begin{pmatrix}
    V\\ V_{'}
  \end{pmatrix}=
  \left(
    \begin{pmatrix}\mathbf{e}\\ \mathbf{0}\end{pmatrix},\;\;
    \begin{pmatrix}X & O\\ I & X\end{pmatrix}
    \begin{pmatrix}\mathbf{e}\\ \mathbf{0}\end{pmatrix},\;\;
    \begin{pmatrix}X & O\\ I & X\end{pmatrix}^2
    \begin{pmatrix}\mathbf{e}\\ \mathbf{0}\end{pmatrix},\;\;\ldots,
    \begin{pmatrix}X & O\\ I & X\end{pmatrix}^n
    \begin{pmatrix}\mathbf{e}\\ \mathbf{0}\end{pmatrix}
  \right).
\]
So, the Arnoldi process applies and can be used to generate a matrix $Q$ with orthogonal columns and
a Hessenberg matrix $H$. Then, the evaluation matrix $W$ can be constructed with
\[
  \text{the first column }\begin{pmatrix}\mathbf{e}\\\mathbf{0}\end{pmatrix},\quad
  \begin{pmatrix}S & O\\I & S\end{pmatrix} W_{-} = WH.
\]
Following \cite{va2021}, we give the MATLAB codes for the method of confluent Vandermonde with
Arnoldi.
\begin{verbatim}
     function [d,H] = polyfitAh(x,f,fp,n)
     m = size(x,1); Q = zeros(2*m,n+1); Q(1:m,1) = 1; H = zeros(n+1,n);
     for k = 1:n
         q = [x.*Q(1:m,k); Q(1:m,k)+x.*Q(m+1:2*m,k)];
         for j = 1:k
             H(j,k) = Q(:,j)'*q/m; q = q - H(j,k)*Q(:,j);        
         end
         H(k+1,k) = norm(q)/sqrt(m); Q(:,k+1) = q/H(k+1,k);
     end
     d = Q\[f;fp];

     function [y,yp] = polyvalAh(d, H, s)
     M = size(s,1);  n = size(d,1)-1; W = zeros(2*M,n+1); W(1:M,1) = 1;
     for k = 1:n
         w = [s.*W(1:M,k); W(1:M,k)+s.*W(M+1:2*M,k)];
         for j = 1:k
             w = w - H(j,k)*W(:,j);
         end
         W(:,k+1) = w/H(k+1,k);
     end
     y = W*d;  yp = y(M+1:2*M,:);  y = y(1:M,:);

\end{verbatim}
To compare, we give also the naive method of confluent Vandermonde without Arnoldi.
\begin{verbatim}
     function c = polyfith(x,f,fp,n)
     A = x.^(0:n); o = zeros(size(x)); D = diag(1:n); A = [A; o A(:,1:n)*D]; c = A\[f;fp];

     function [y,yp] = polyvalh(c,s)
     n = length(c)-1;  B = s.^(0:n);  y = B*c;
     o = zeros(size(s));  D = diag(1:n);  yp = [o B(:,1:n)*D]*c;
\end{verbatim}

Other than the Hermite scenario for both fitting and evaluation, we can also fit only the function
values using \verb|polyfitA| (or \verb|polyfit|) from \cite{va2021} but evaluates the function
values and derivatives using \verb|polyvalAh| (or \verb|polyvalh|) here. This scenario is more
interesting if the data of derivatives are not available for fitting. In other words, we construct
first a discrete orthogonal basis for the function values by \verb|polyfitA| and then the derivative
matrix \verb|W(M+1:2*M,:)| under the same basis by \verb|polyvalAh|. We notice that an equivalent
formula for doing this has been given in \cite{lt2022}, which we should detail in the Discussion
section. There is a third possibility: one can still construct a discrete orthogonal basis for the
Hermite data, but fit only the function values, and then construct the derivative evaluation matrix
under that basis; however we found no advantage of it. The fourth possibility is to use the function
values orthogonalized basis for Hermite fitting and evaluation, which we found performs equally well
as using the Hermite data orthogonalized basis for Examples 1--2 in the following sections. We will
see also for the other examples the comparison results of Hermite and function values orthogonalized
bases.

Both the fitting and evaluation can be generalized to high-order derivatives. For example, let
$V_{''}$ be the second-order derivative matrix under the standard basis. We have
\[
  \begin{pmatrix}
    V\\ V_{'}\\ V_{''}
  \end{pmatrix}=
  \left(
    \begin{pmatrix}\mathbf{e}\\ \mathbf{0}\\ \mathbf{0}\end{pmatrix},\;\;
    \begin{pmatrix}X & O & O \\ I & X & O\\ O & 2I & X\end{pmatrix}
    \begin{pmatrix}\mathbf{e}\\ \mathbf{0}\\ \mathbf{0}\end{pmatrix},\;\;
    \begin{pmatrix}X & O & O\\ I & X & O\\ O & 2I & X\end{pmatrix}^2
    \begin{pmatrix}\mathbf{e}\\ \mathbf{0}\\ \mathbf{0}\end{pmatrix},\;\;\ldots,
    \begin{pmatrix}X & O & O\\ I & X & O\\ O & 2I & X\end{pmatrix}^n
    \begin{pmatrix}\mathbf{e}\\ \mathbf{0}\\ \mathbf{0}\end{pmatrix}
  \right).
\]
It may be worthwhile to state the general case as follows.
\begin{theorem}
  Let $V$ be the Vandermonde matrix and $V_{(\ell)}$ be the $\ell$th order derivative matrix
  corresponding to $V$. Then, the matrix $[V; V_{(1)}; V_{(2)}; \ldots; V_{(\ell)}]$ has the $(k+1)$th
  column exactly equal to
\[
  \begin{pmatrix}
    X & O & \cdots & \cdots & O \\
    I & X & \ddots & \vdots & \vdots\\
    O & 2I & \ddots & \ddots & \vdots\\
    \vdots & \ddots & \ddots & \ddots & O\\
    O & \cdots & O & \ell I & X
  \end{pmatrix}^{k}
  \begin{pmatrix}\mathbf{e}\vphantom{X}\\ \mathbf{0}\vphantom{\ddots}\\
    \vdots\vphantom{X}\\\mathbf{0}\vphantom{\ddots}\\ \mathbf{0}\vphantom{X}\end{pmatrix}.
\]  
\end{theorem}

\section{Example 1: interpolation in Chebyshev points.}

\begin{figure}
  \centering
  \includegraphics[scale=.16]{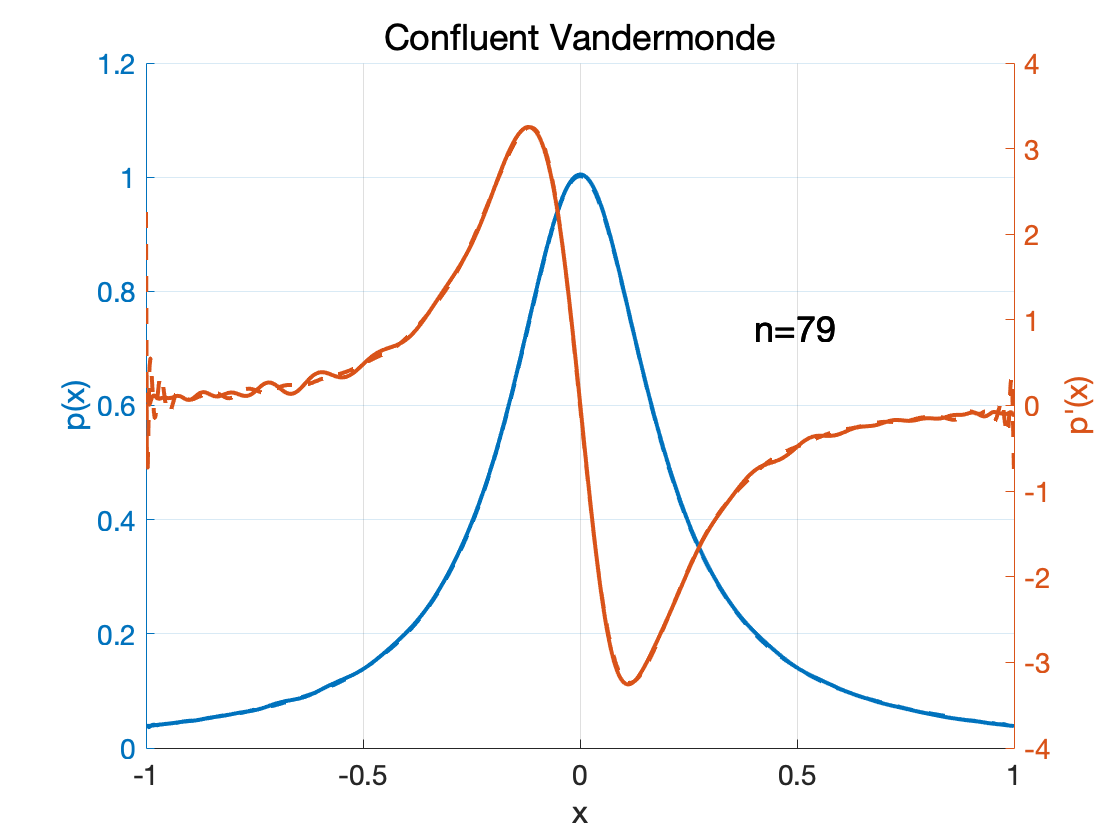}\quad%
  \includegraphics[scale=.16]{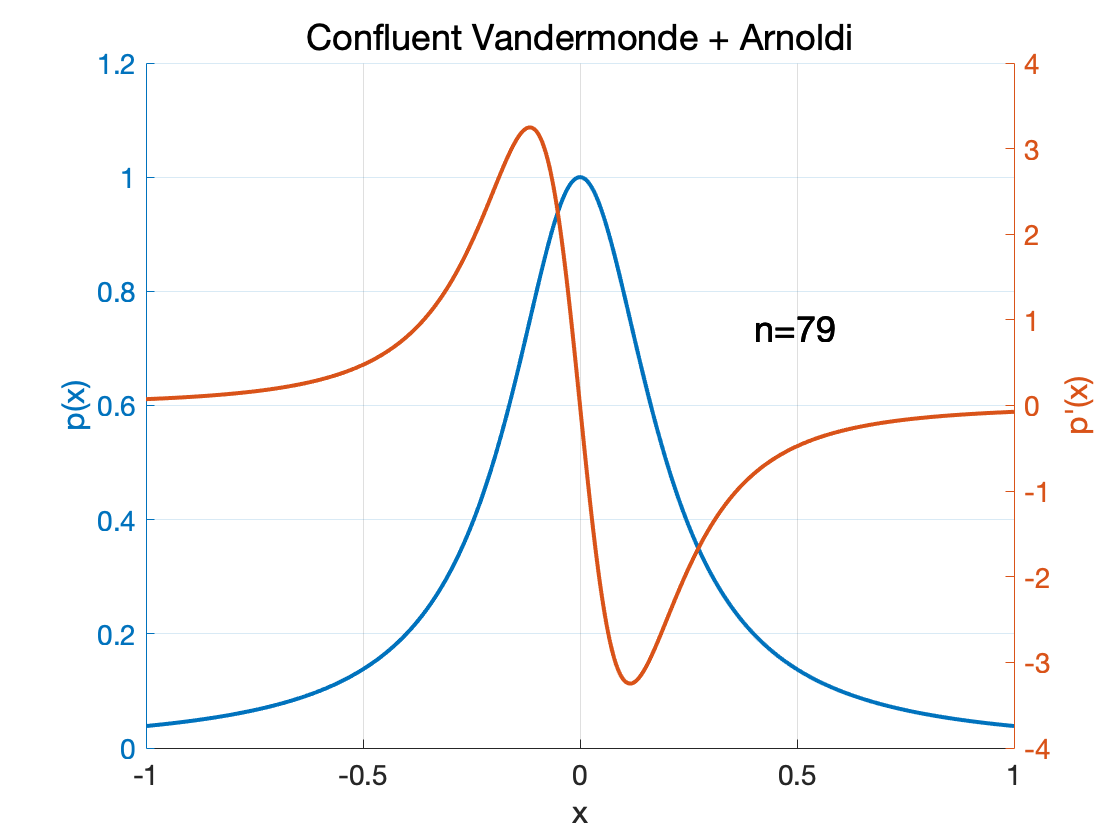}\\
  \includegraphics[scale=.16]{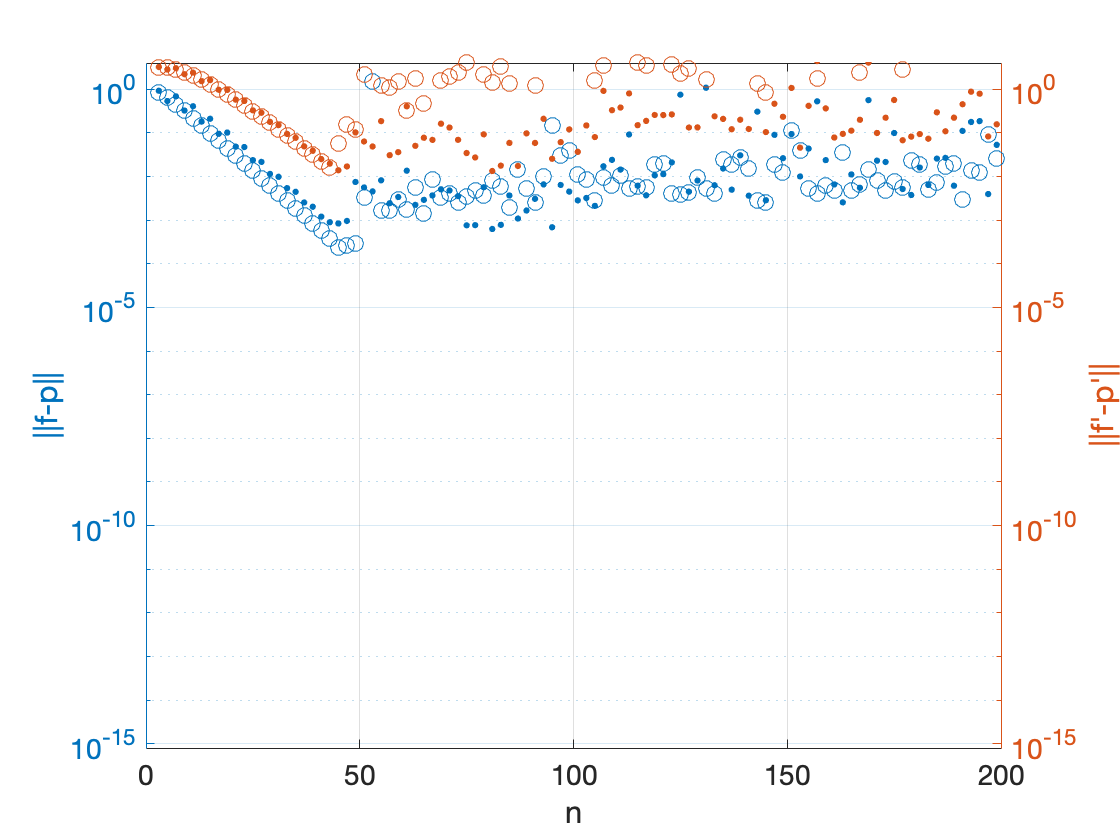}\quad%
  \includegraphics[scale=.16]{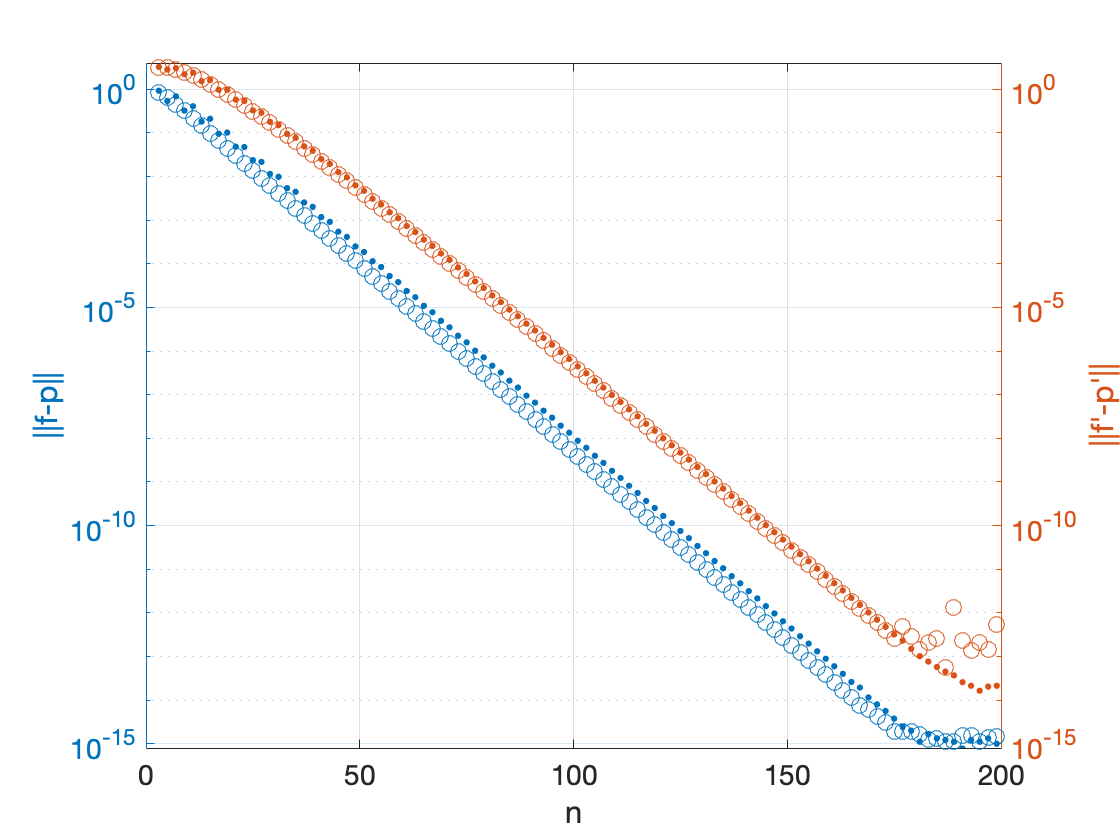}
  \caption{\small{Hermite (solid lines/points) and Lagrange (dashed lines/circles) interpolation of
    $f(x)=1/(1+25x^2)$ on Chebyshev grids of $[-1,1]$ by the confluent Vandermonde without (left
    column) and with (right column) Arnoldi.}}
  \label{figex1}
\end{figure}

This example is used in \cite{va2021} for the Lagrange interpolation. Now we interpolate the
function $f(x)=1/(1+25x^2)$ and its derivative simultaneously (\verb|polyfitAh|+\verb|polyvalAh|) by
a degree $n$ (odd) polynomial in $m=(n+1)/2$ Chebyshev points $x_j=\cos((m-j)\pi/(m-1))$,
$1\le j\le m$. We also do the Lagrange interpolation in $n+1$ points as in \cite{va2021} and
evaluate the derivative as proposed here (\verb|polyfitA|+\verb|polyvalAh|). The results are shown
in Figure~\ref{figex1}. We see that with Arnoldi both $f$ and $f'$ can be approximated to high
accuracy with the exponential convergence rate, while without Arnoldi the method is prone to
round-off errors as $n>50$. The derivative evaluation enabled by the confluent Vandermonde with
Arnoldi (\verb|polyvalAh|) works equally well for the Hermite and Lagrange interpolation. Throughout
the paper, the function norm is the maximum norm computed on a very fine grid.

\section{Example 2: least-squares fitting on two intervals.}

\begin{figure}
  \centering
  \includegraphics[scale=.16]{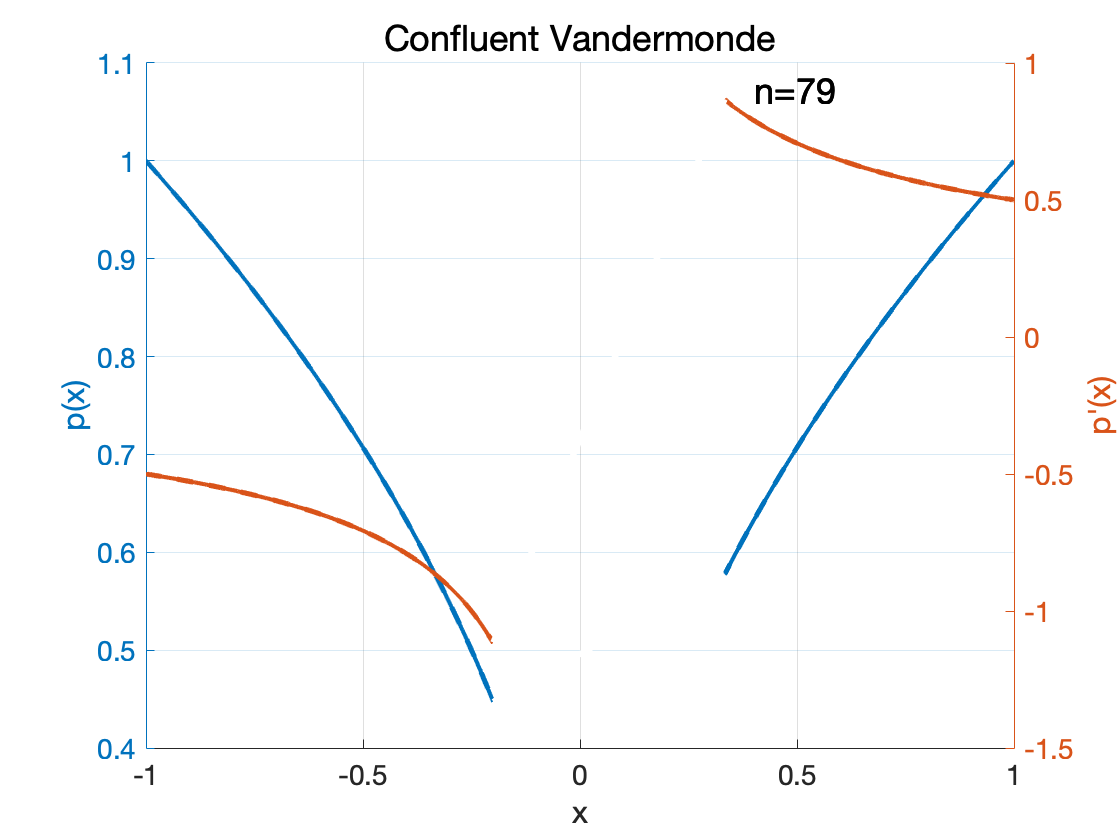}\quad%
  \includegraphics[scale=.16]{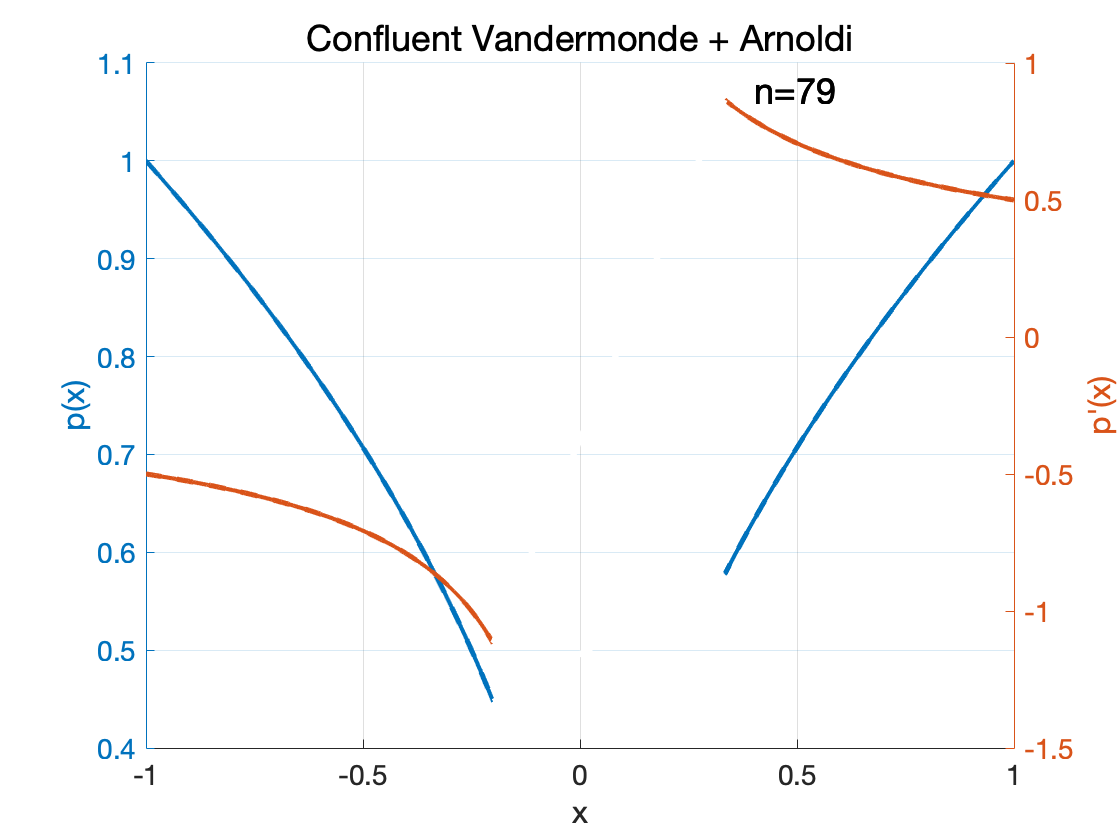}\\
  \includegraphics[scale=.16]{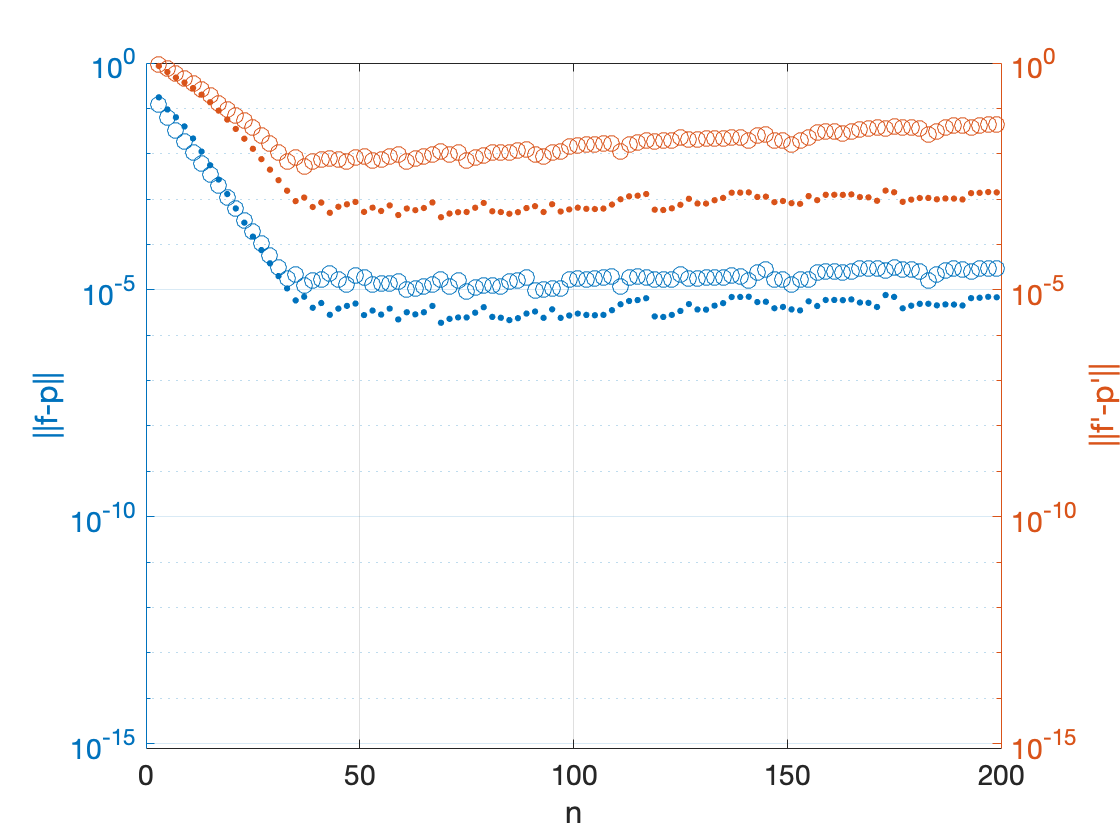}\quad%
  \includegraphics[scale=.16]{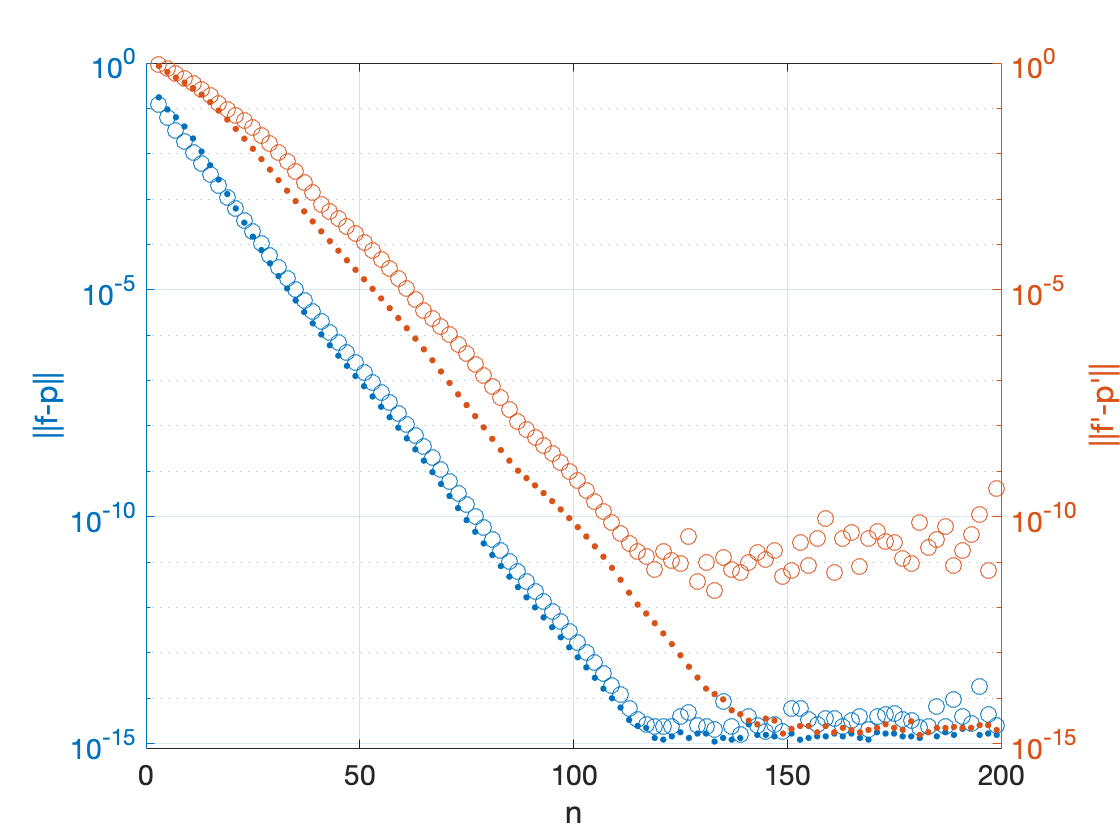}
  \caption{\small{Least-squares fitting of $f(x)=\sqrt{|x|}$ on $[-1,-1/3]\cup [1/5,1]$ (dashed
      lines/circles) and simultaneously its derivatives (solid lines/points) in equally spaced
      points by the confluent Vandermonde without (left column) and with (right column) Arnoldi.}}
  \label{figex2}
\end{figure}

For the function $f(x)=\sqrt{|x|}$ on $[-1,-1/3]\cup [1/5,1]$, we sample the function values (and
its derivatives for the Hermite version) on every interval with equally spaced points, and fit by
least-squares these data with a degree $n$ polynomials. The number of sampling points $m$ on each
interval is chosen as $m=5(n+1)$ for the Hermite version and $m=10(n+1)$ for the Lagrange
version. The results are shown in Figure~\ref{figex2}.

\section{Example 3: Hermite Fourier extension.}

\begin{figure}
  \centering
  \includegraphics[scale=.16]{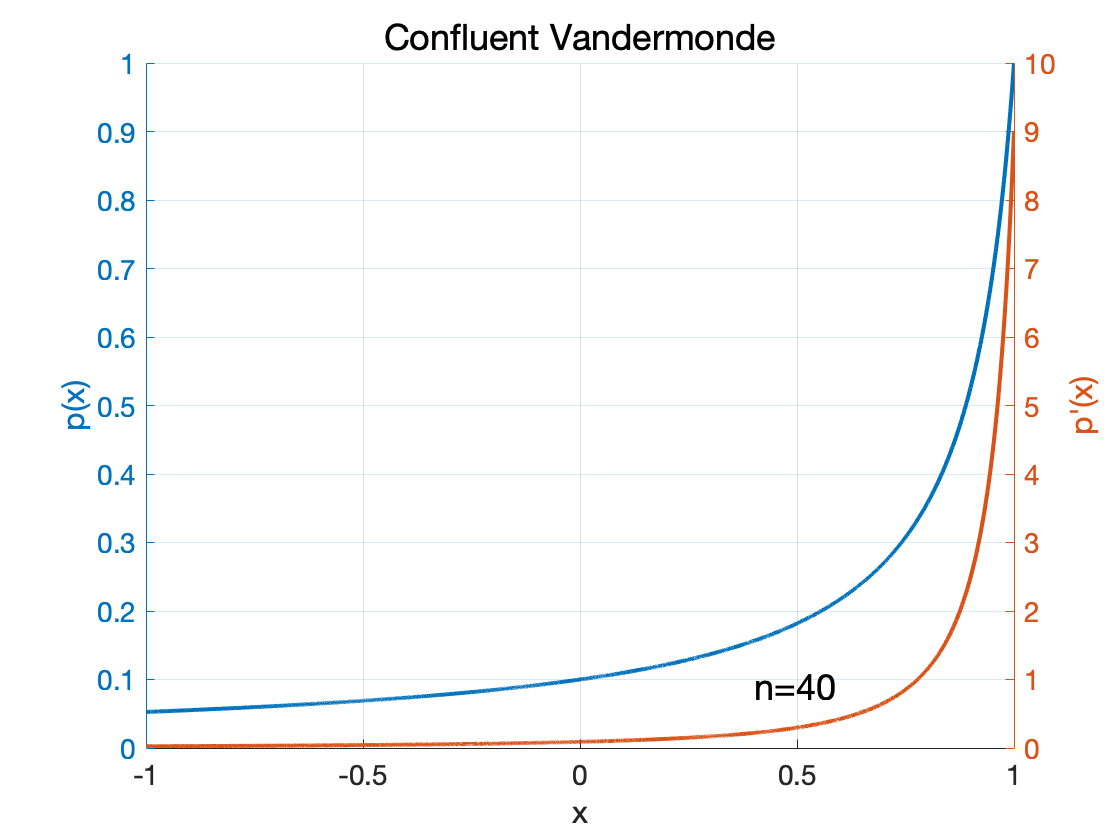}\quad%
  \includegraphics[scale=.16]{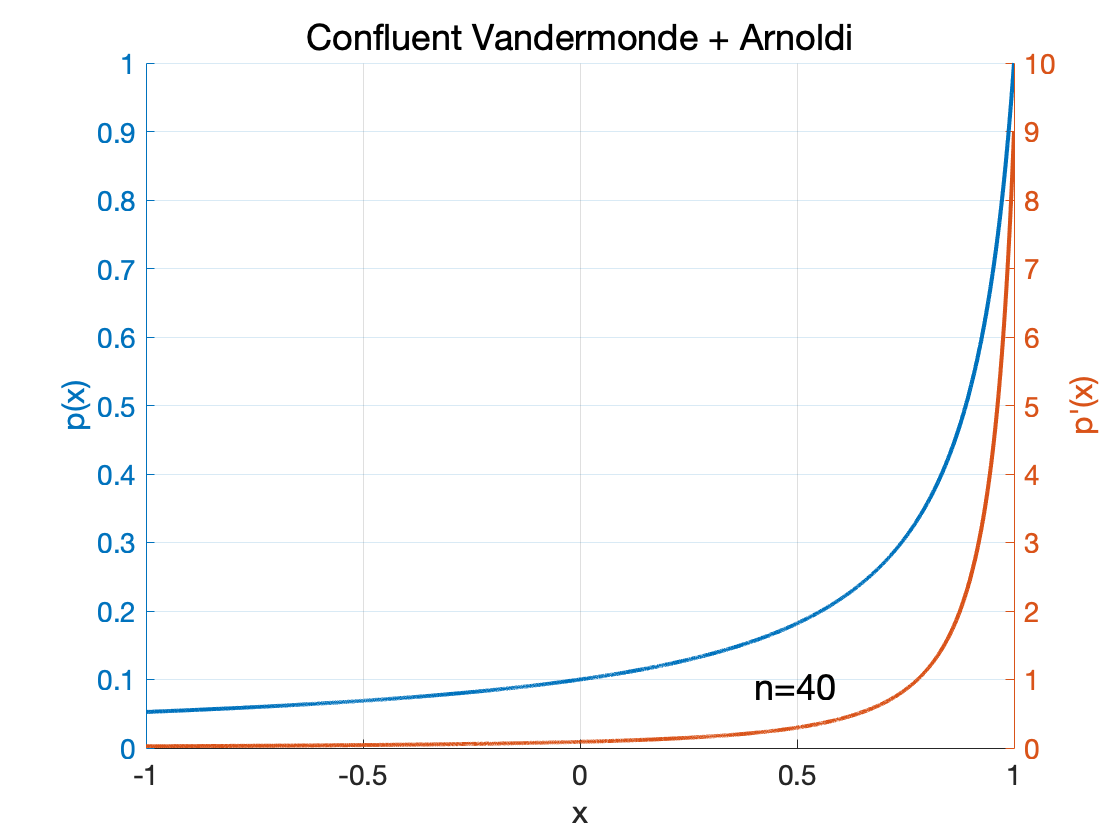}\\
  \includegraphics[scale=.16]{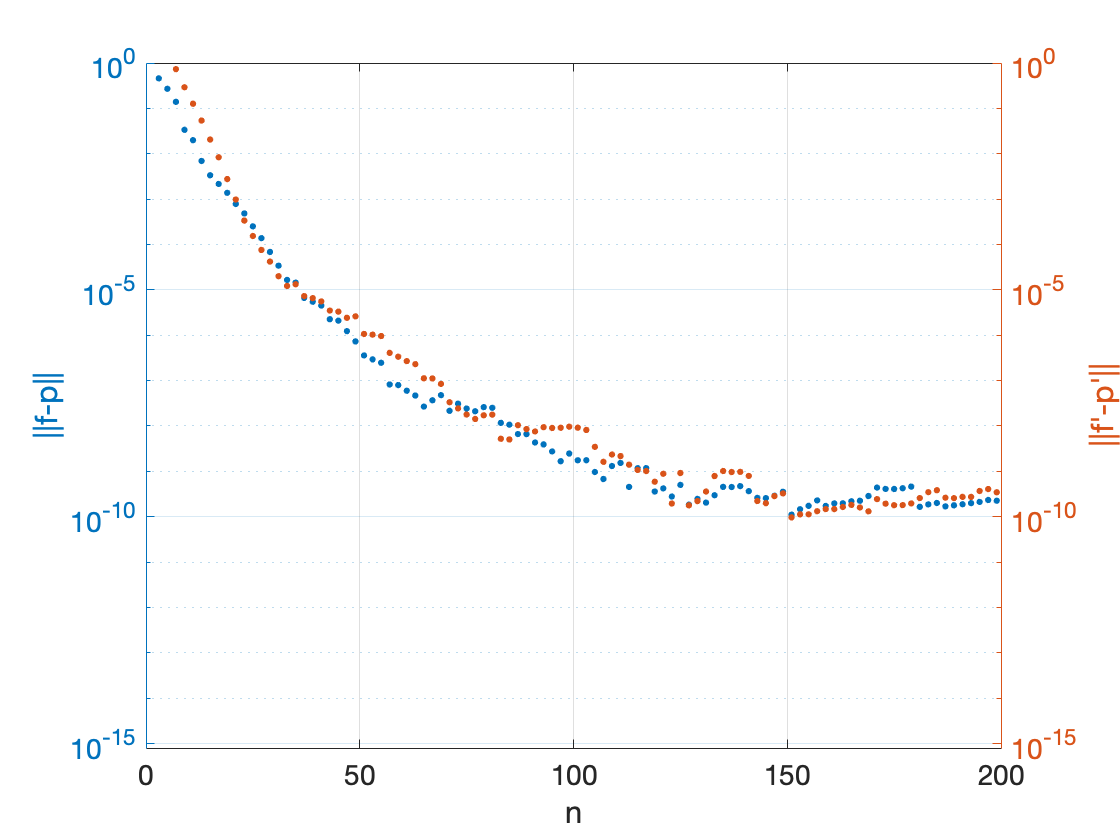}\quad%
  \includegraphics[scale=.16]{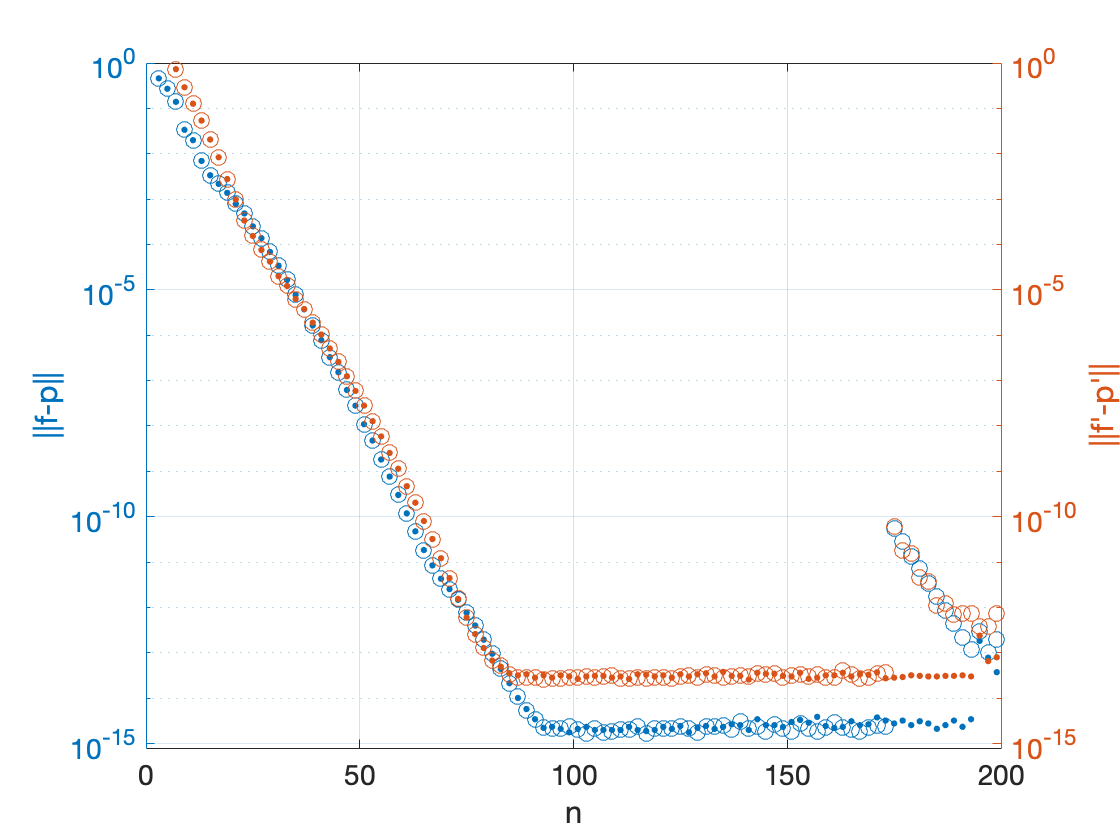}
  \caption{\small{Fourier extension of
      $f(x)=1/(10-9x)$ and simultaneously its derivatives on Chebyshev grids of
      $[-1,1]$ by the truncated Fourier series on
      $[-2,2]$ using the confluent Vandermonde without (left column) and with (right column)
      Arnoldi. Bottom right: Hermite (dots) vs function values (circles) orthogonalized bases.}}
  \label{figex3}
\end{figure}

This example is used in \cite{va2021} for approximation of $f(x)=1/(10-9x)$ on $[-1,1]$ by the
Fourier series on the larger interval $[-2,2]$. Now we consider both the function values and its
derivatives:
\[
  f(x)\approx\mathrm{Re}\left(\sum_{k=0}^n\mathrm{e}^{\mathrm{i}k\pi x/2}c_k\right),\quad
  -\frac{2}{\pi}f'(x)\approx\mathrm{Im}\left(\sum_{k=0}^nk\mathrm{e}^{\mathrm{i}k\pi x/2}c_k\right).
\]
Following \cite{va2021}, this is realized with $z:=\mathrm{e}^{\mathrm{i}\pi x/2}$,
$c_k=a_k+\mathrm{i}b_k$ ($a_k, b_k\in\mathbb{R}$) and
\[
  f(x) \approx\sum_{k=0}^n(\mathrm{Re}z^ka_k-\mathrm{Im}z^kb_k),\quad
  -\frac{2}{\pi}f'(x) \approx\sum_{k=0}^n\left(\mathrm{Im}(kz^k)a_k+\mathrm{Re}(kz^k)b_k\right).
\]
The last line of \verb|polyfith| is changed to
\begin{verbatim}
     m = length(x);  c = [real(A(1:m,:))        imag(A(1:m,2:n+1)); ...
                          imag(A(m+1:2*m,:))   -real(A(m+1:2*m,2:n+1))]\[f;fp];
     c = c(1:n+1) - 1i*[0; c(n+2:2*n+1)];
\end{verbatim}
where \verb|fp| stores the scaled derivatives
$-\frac{2}{\pi}f'(x_j)$. Similar changes should be made also in \verb|polyfitAh| but with
$Q$ instead of
$A$. In \verb|polyvalh| and \verb|polyvalAh|, now we take the real part for \verb|y| and imaginary
part for \verb|yp|. Note again that the output \verb|yp| stores the scaled derivatives
$-\frac{2}{\pi}p'(s_j)$. Following \cite{va2021}, 500 Chebyshev points on
$[-1,1]$ are used for the Hermite fitting. The results are shown in Figure~\ref{figex3}. The
Lagrange fitting is not present in the figure since the derivative evaluation for it gave totally
wrong results.

\section{Example 4: Laplace equation and Dirichlet-to-Neumann map.}

\begin{figure}
  \centering
  \includegraphics[scale=.16]{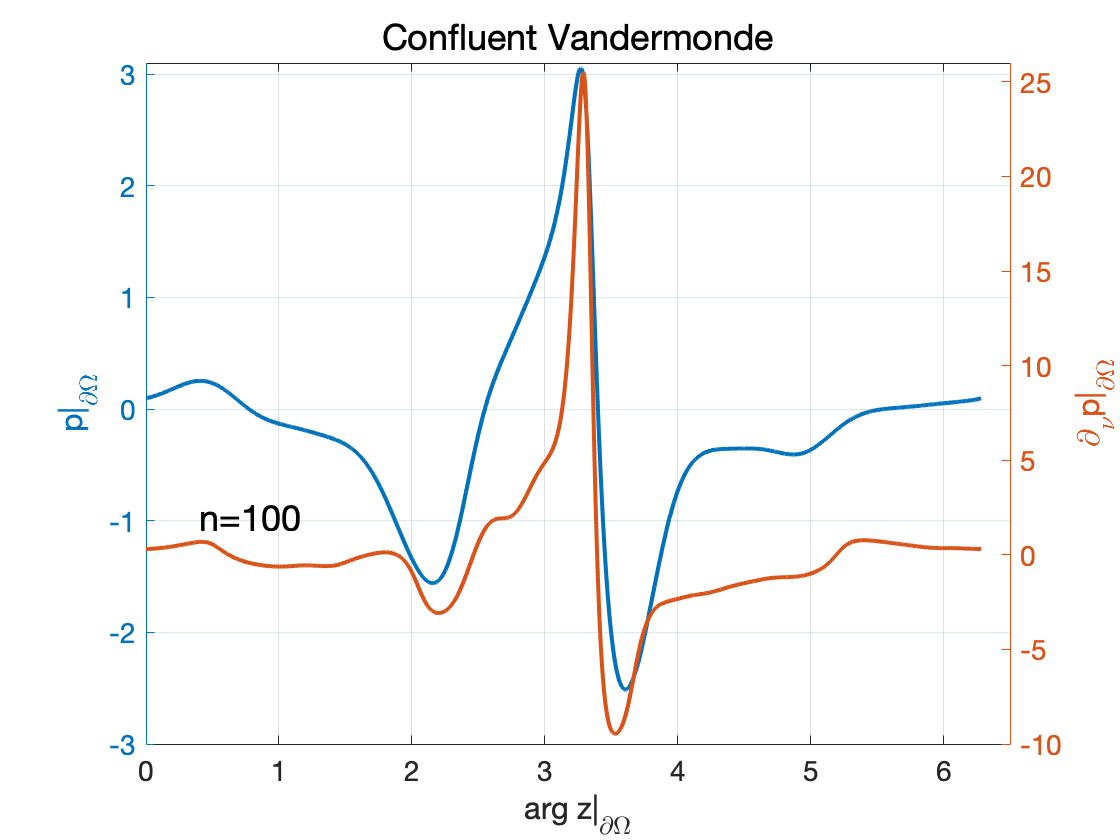}\quad%
  \includegraphics[scale=.16]{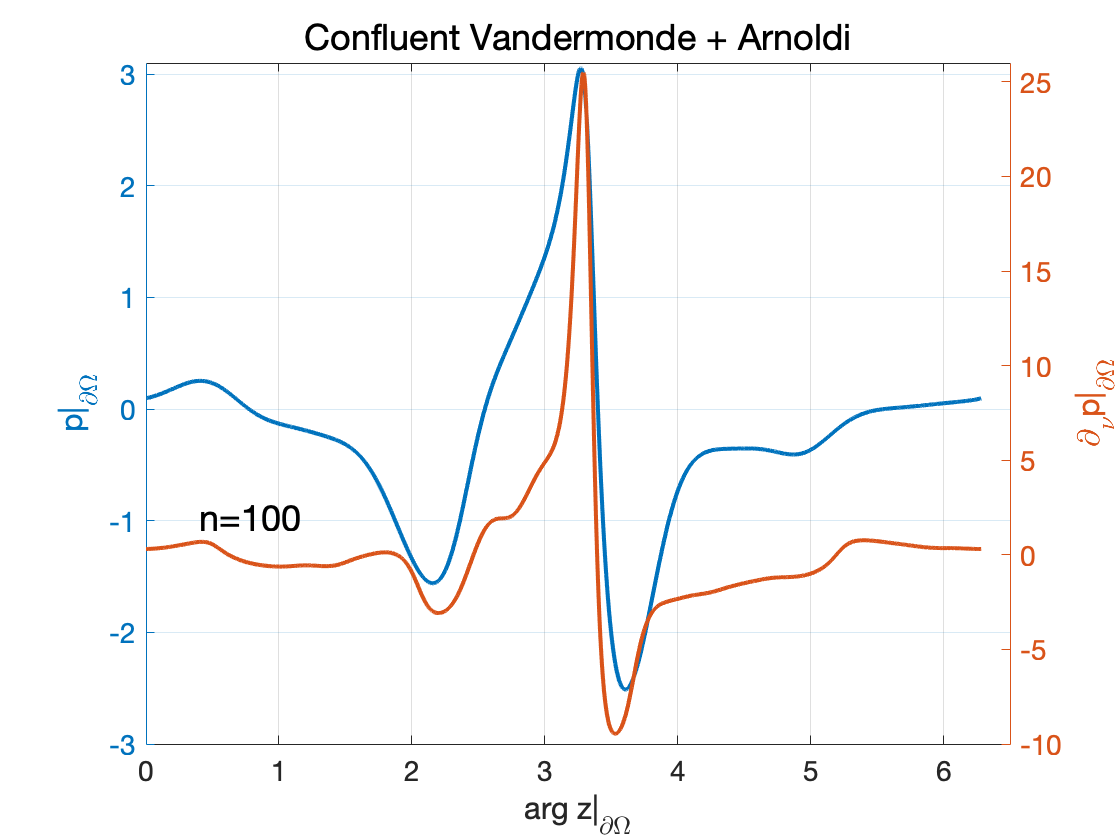}\\
  \includegraphics[scale=.16]{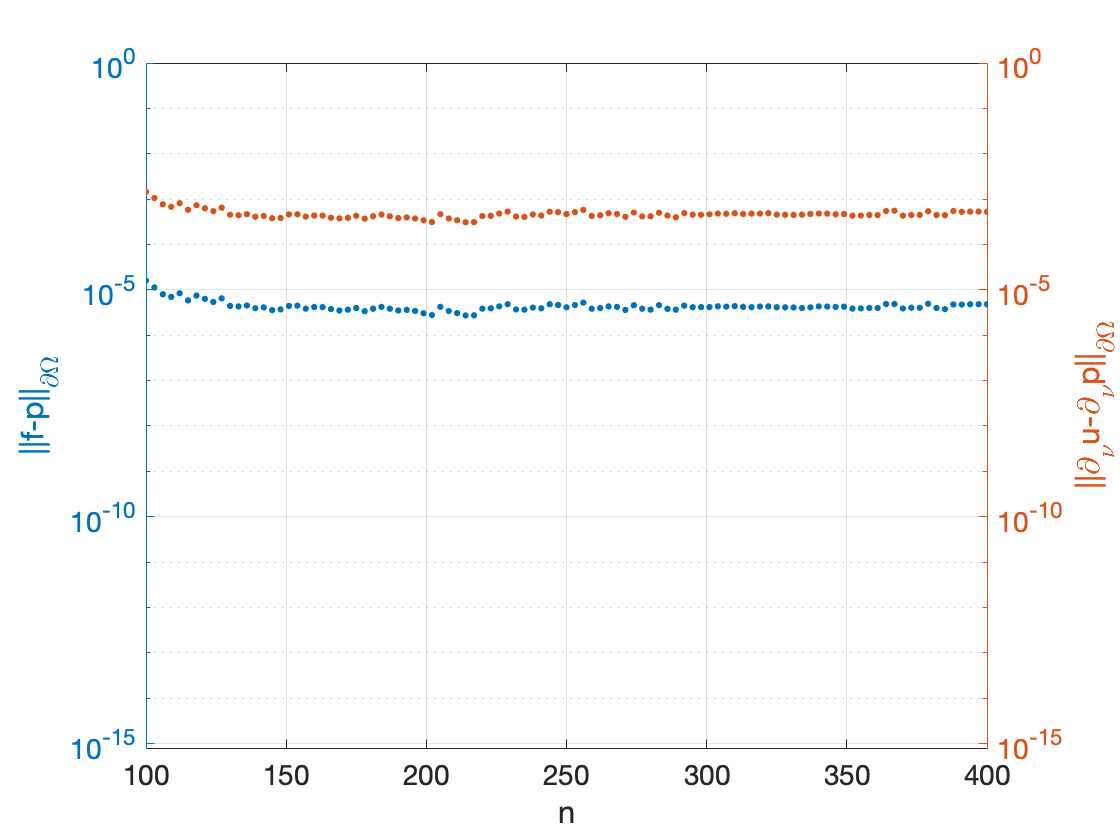}\quad%
  \includegraphics[scale=.16]{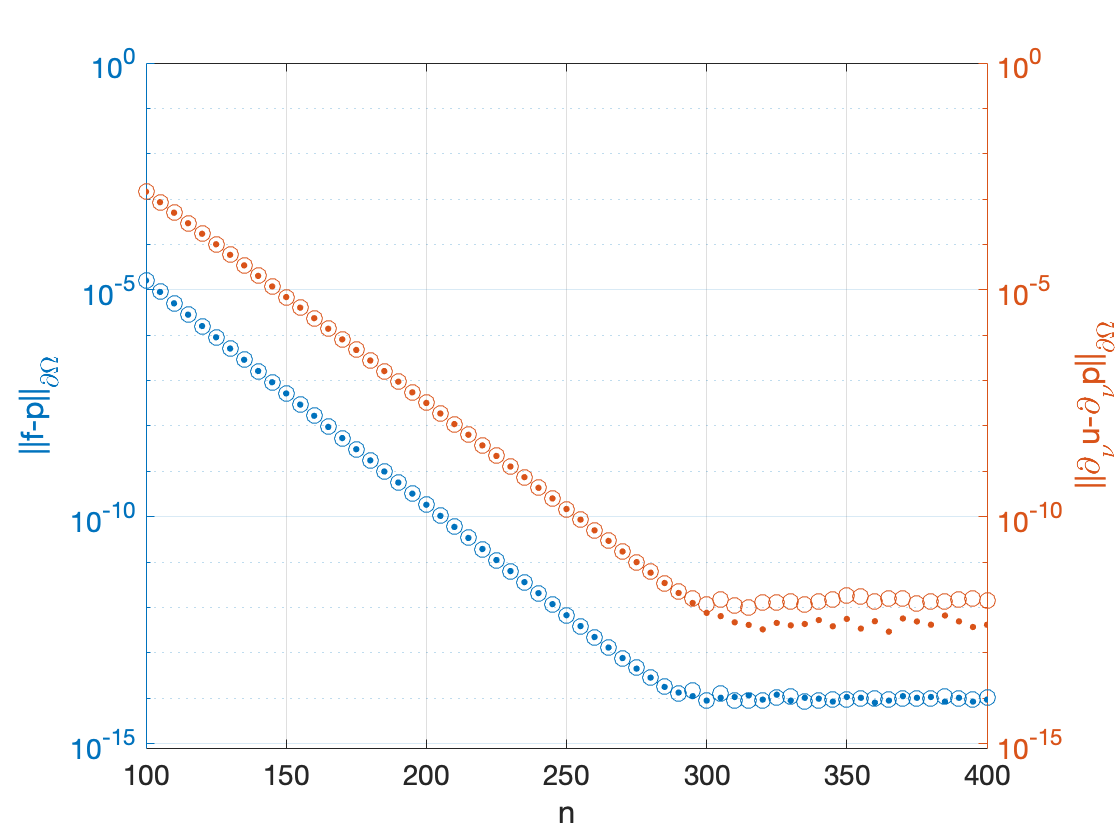}
  \caption{\small{Dirichlet-to-Neumann map of the Dirichlet data $f(z)=\mathrm{Re}\,(\log(0.8+z))^2$
      at
      $z\in\partial\Omega=\{\mathrm{e}^{\mathrm{i}t}(0.7+0.25\cos(4t-2)+0.05\cos(8t-4)): 0\le t=\arg
      z\le 2\pi\}$ by the confluent Vandermonde without (left column) and with (right column)
      Arnoldi. Bottom-right: Hermite (dots) vs function values (circles) orthogonalized bases.}}
  \label{figex4}
\end{figure}

Given a domain $\Omega\subset\mathbb{R}^2$ with boundary $\partial\Omega$ and a function
$f:\partial\Omega\to \mathbb{R}$, the Dirichlet-to-Neumann map $T$ involves first solving the
Laplace equation $-\Delta u =0$ in $\Omega$, $u=f$ on $\partial\Omega$ and then evaluating the outer
normal derivative $\partial_{\boldsymbol\nu}u=:Tf$ on $\partial\Omega$. By identifying
$(x,y)\in\mathbb{R}^2$ and $z=x+\mathrm{i}y\in\mathbb{C}$, $u$ can be viewed as the real part of an
analytic function $h$ in the complex plane, so that $-\Delta u=0$ is automatically
satisfied. Following \cite{va2021}, we approximate $h(z)$ by a polynomial $\sum_{k=0}^nc_kz^k$. Let
$c_k=a_k+\mathrm{i}b_k$ ($a_k$, $b_k$ $\in\mathbb{R}$). Then,
$u(z) \approx\sum_{k=0}^n(\mathrm{Re}z^ka_k-\mathrm{Im}z^kb_k)$. To compute $Tf$, we need only to
fit the boundary condition $u=f$ on $\partial\Omega$ and then evaluate
$\partial_{\boldsymbol\nu}u=\mathrm{Re}(\nu(z)h'(z))$ (since
$h'(z)=\partial_xu-\mathrm{i}\partial_yu$ and $\nu(z)=\nu_1(z) + \mathrm{i}\nu_2(z)$ corresponds to
the unit outer normal $\boldsymbol{\nu}=(\nu_1,\nu_2)$ on $\partial\Omega$), or more precisely
\[\partial_{\boldsymbol{\nu}}u=
  \sum_{k=1}^n(\nu_1(z)\mathrm{Re}z^{k-1}k-\nu_2(z)\mathrm{Im}z^{k-1}k)a_k-
  \sum_{k=1}^n(\nu_2(z)\mathrm{Re}z^{k-1}k+\nu_1(z)\mathrm{Im}z^{k-1}k)b_k.
\]
A concrete example is demonstrated in Figure~\ref{figex4}, where $m=10n$ points with the argument of
$z\in\partial\Omega$ equally spaced are used for the fitting. We see that the Hermite orthogonalized
basis leads to a bit more accurate results than the function values orthogonalized basis does.

\section{Example 5: Steklov eigenvalue problem in an ellipse.}

\begin{figure}
  \centering
  \includegraphics[scale=.16]{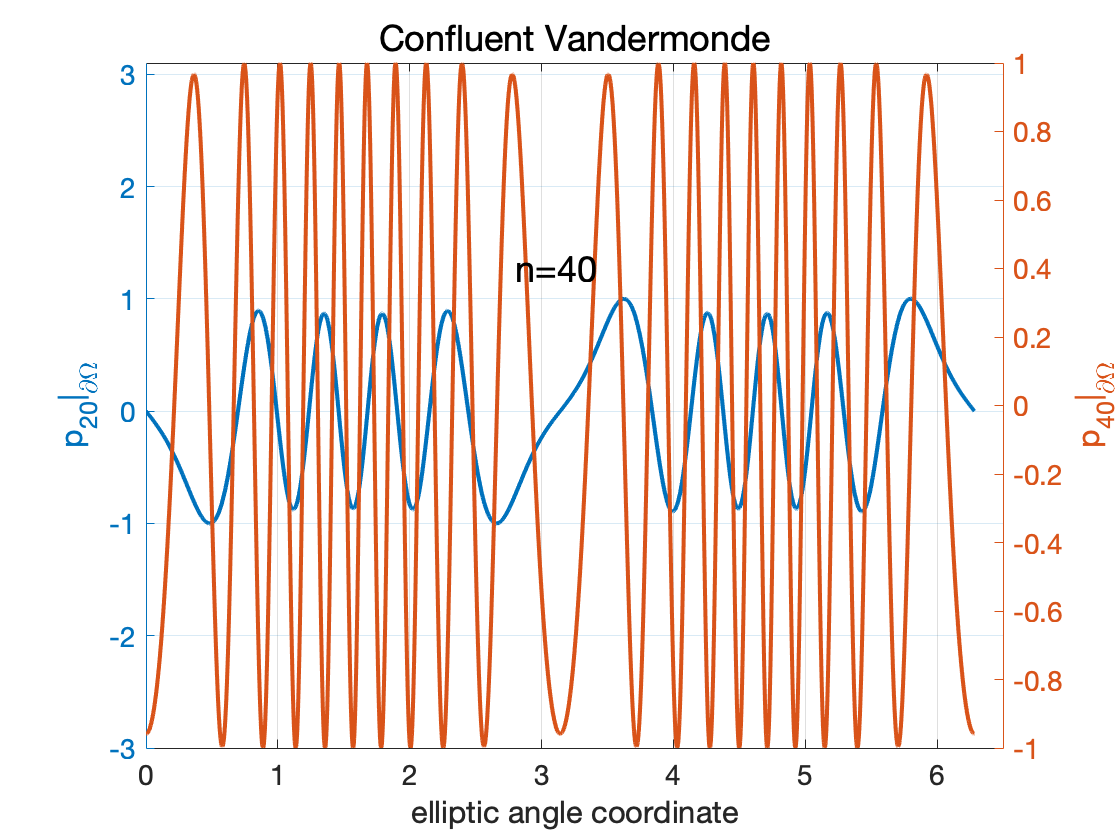}\quad%
  \includegraphics[scale=.16]{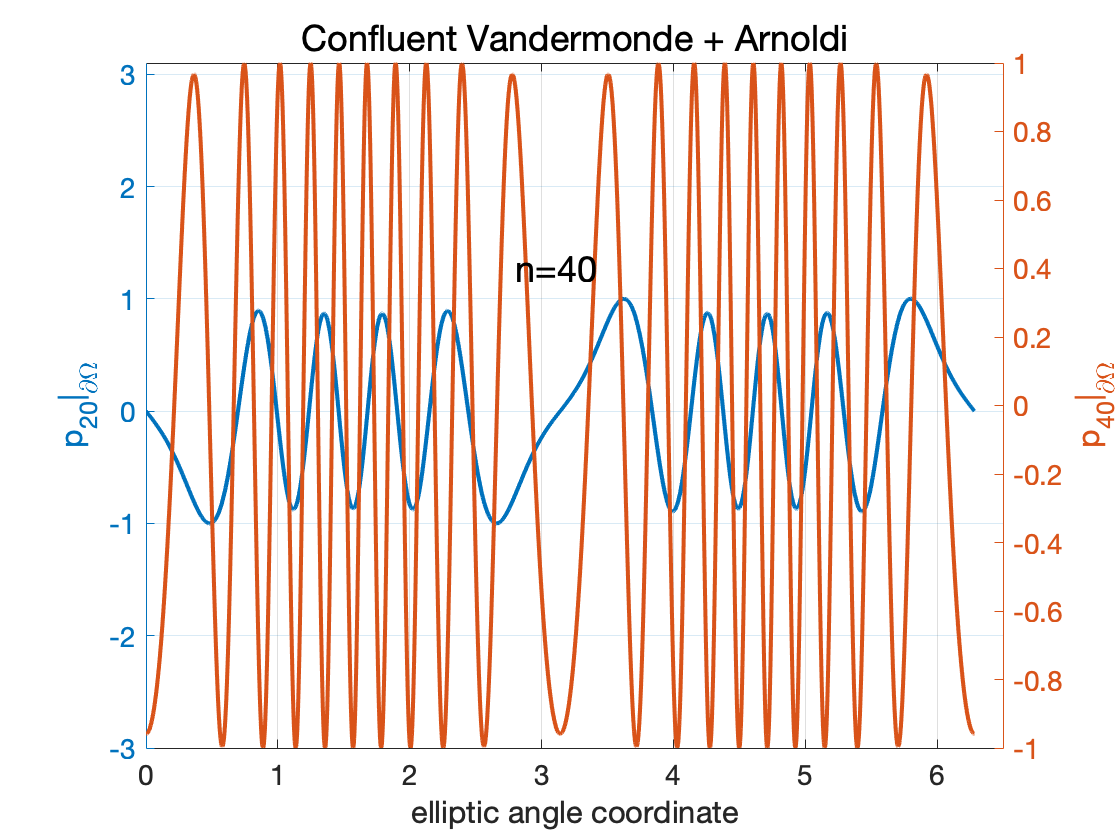}\\
  \includegraphics[scale=.16]{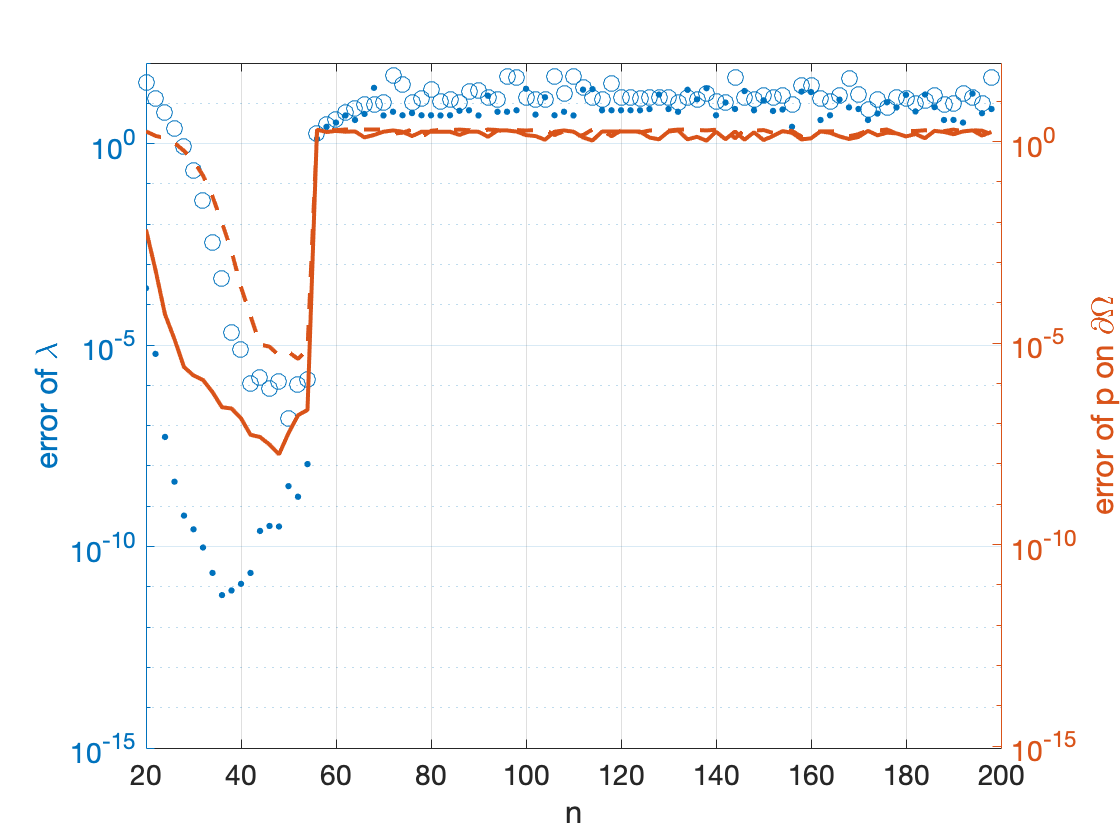}\quad%
  \includegraphics[scale=.16]{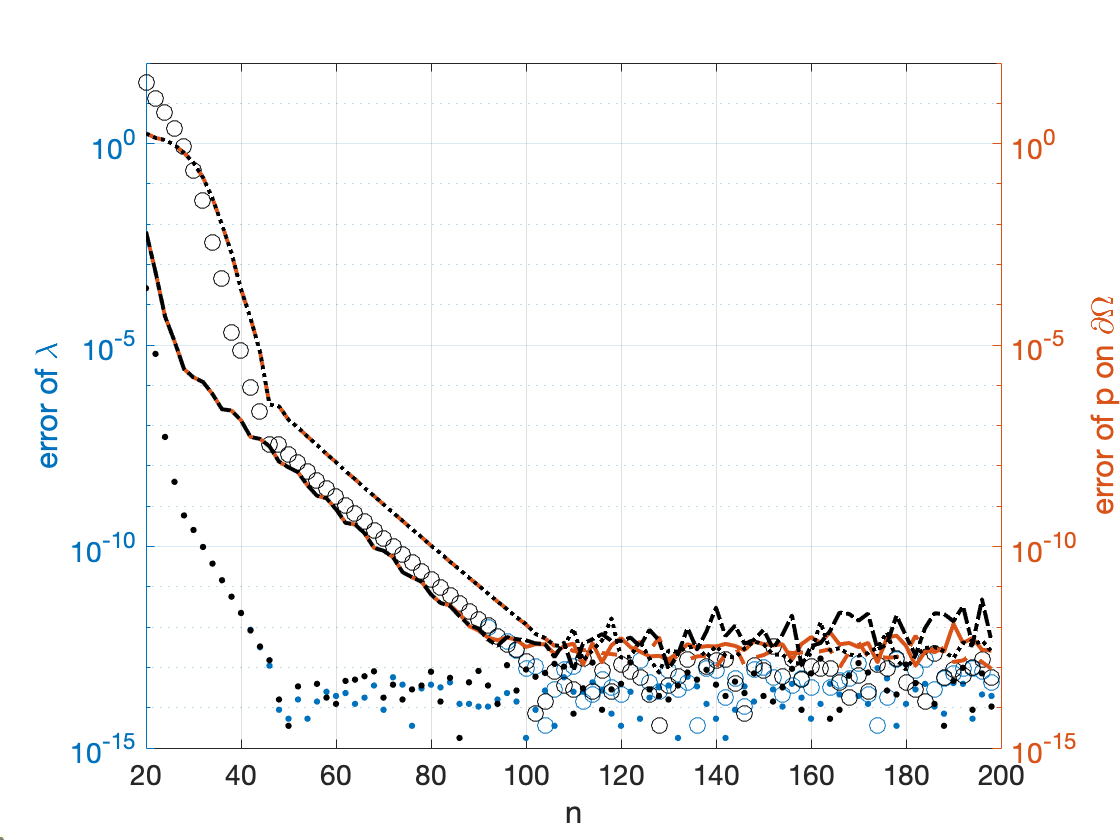}
  \caption{\small{Laplace Steklov eigenvalue problem in the ellipse
      $z=\cos t + \mathrm{i}\frac{1}{5}\sin t$ by the confluent Vandermonde without (left column)
      and with (right column) Arnoldi. Top: the $k$th eigenfunction $p_k$ ($\|p_k\|=1$), $k=20, 40$,
      for the eigenvalues $0=\lambda_1\le \ldots\le \lambda_k\le\ldots$. Bottom: errors of $p_{20}$,
      $\lambda_{20}$ ($-$, $\cdot$) and $p_{40}$, $\lambda_{40}$ ($--$, $\circ$). Bottom right:
      Hermite vs function values (in black) orthogonalized bases.}}
  \label{figex5}
\end{figure}

This is a continuation of the previous section. We look for the eigenvalue $\lambda$ of $T$ and the
assoicated eigenfunction $f:\partial\Omega\to \mathbb{R}$ such that $Tf=\lambda f$. In other words,
we want the real part $u$ of an analytic function in $\Omega$ to satisfy
$\partial_{\boldsymbol\nu}u=\lambda u$ on $\partial\Omega$. Let $Q$ be the matrix from the Arnoldi
process on the confluent Vandermonde matrix; see the Introduction section.  We partition the $2m$
rows of $Q$ in half: $Q=(Q_0; Q_1)$ with $Q_0$ containing the first $m$ rows that represents the
sampling of basis function values and $Q_1$ representing the sampling of derivatives of basis
functions. Let the matrices be separated into real parts and imaginary parts:
$Q_0=Q_{0,r}+\mathrm{i}Q_{0,i}$, $Q_1=Q_{1,r}+\mathrm{i}Q_{1,i}$, and $\nu_1$, $\nu_2$ be diagonal
matrices containing the sampled values of the unit outer normal of $\partial\Omega$. Then, the
discrete Steklov eigenvalue problem is to find $\mathbf{a}\in\mathbb{R}^{n+1}$,
$\mathbf{b}=(0, b_1,\ldots,b_n)^T$ and $\lambda\in\mathbb{R}$ such that
\[
  (\nu_1Q_{1,r} -\nu_2Q_{1,i})\mathbf{a} -(\nu_1Q_{1,i}
  +\nu_2Q_{1,r})\mathbf{b}=\lambda(Q_{0,r}\mathbf{a}-Q_{0,i}\mathbf{b}),
\]
denoted by $Q_{\nu}\boldsymbol{\beta} = \lambda Q_{00}\boldsymbol{\beta}$ with
$\boldsymbol{\beta}=(\mathbf{a}; -\hat{\mathbf{b}})$ ($\hat{\mathbf{b}}=(b_1,\ldots,b_n)^T$), which
is a rectangular eigenvalue problem when $m>2n+1$. Following \cite{re2021}, we first do the economic
QR decomposition $Q_{00}=\tilde{Q}_{00}\tilde{R}_{00}$, then multiply the eigenvalue equation with
$\tilde{Q}_{00}^*$ to get
$\tilde{Q}_{00}^*Q_{\nu}\boldsymbol{\beta}=\lambda \tilde{R}_{00}\boldsymbol{\beta}$, which is a
generalized eigenvalue problem of square matrices and can be solved by the standard method (using
the MATLAB command \verb|eig|). For the solved example, see Figure~\ref{figex5}. The $m=10n+1$
points are sampled at equally spaced elliptic angle coordinates. The benchmark solution uses
$n=400$.

\section{Example 6: sloshing modes on top of a cup.}

\begin{figure}
  \centering
  \includegraphics[scale=.16]{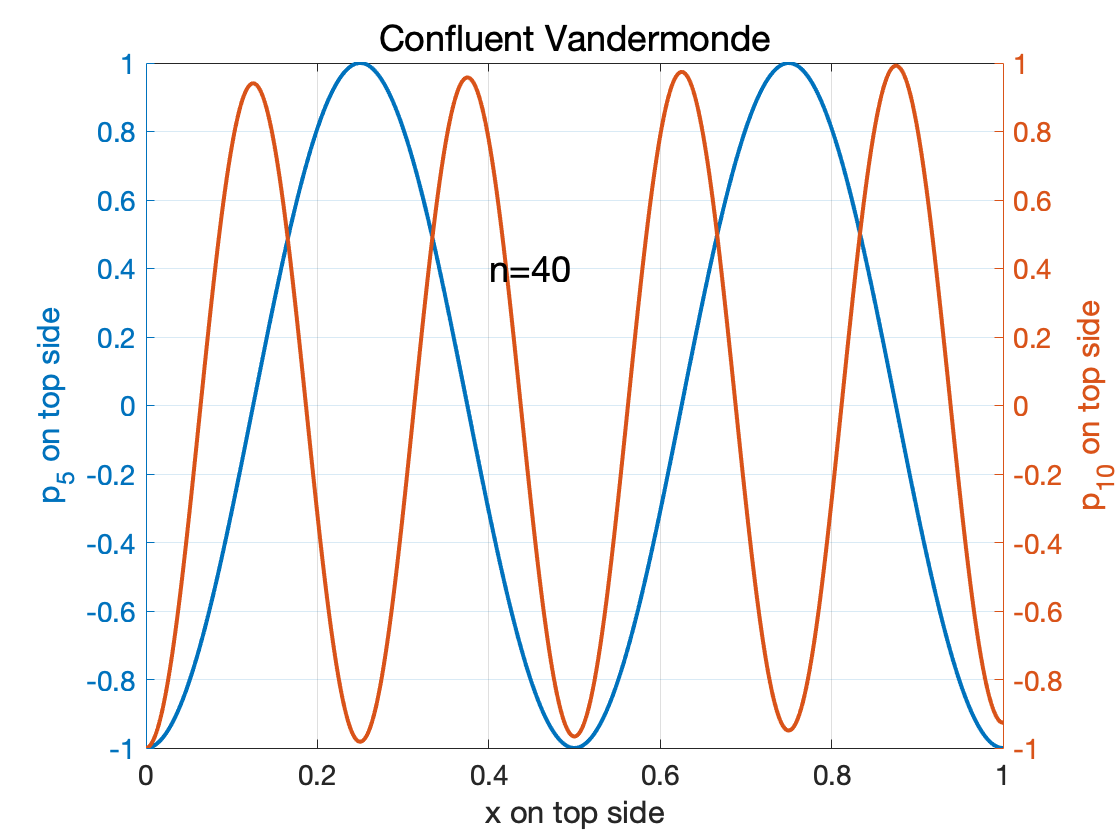}\quad%
  \includegraphics[scale=.16]{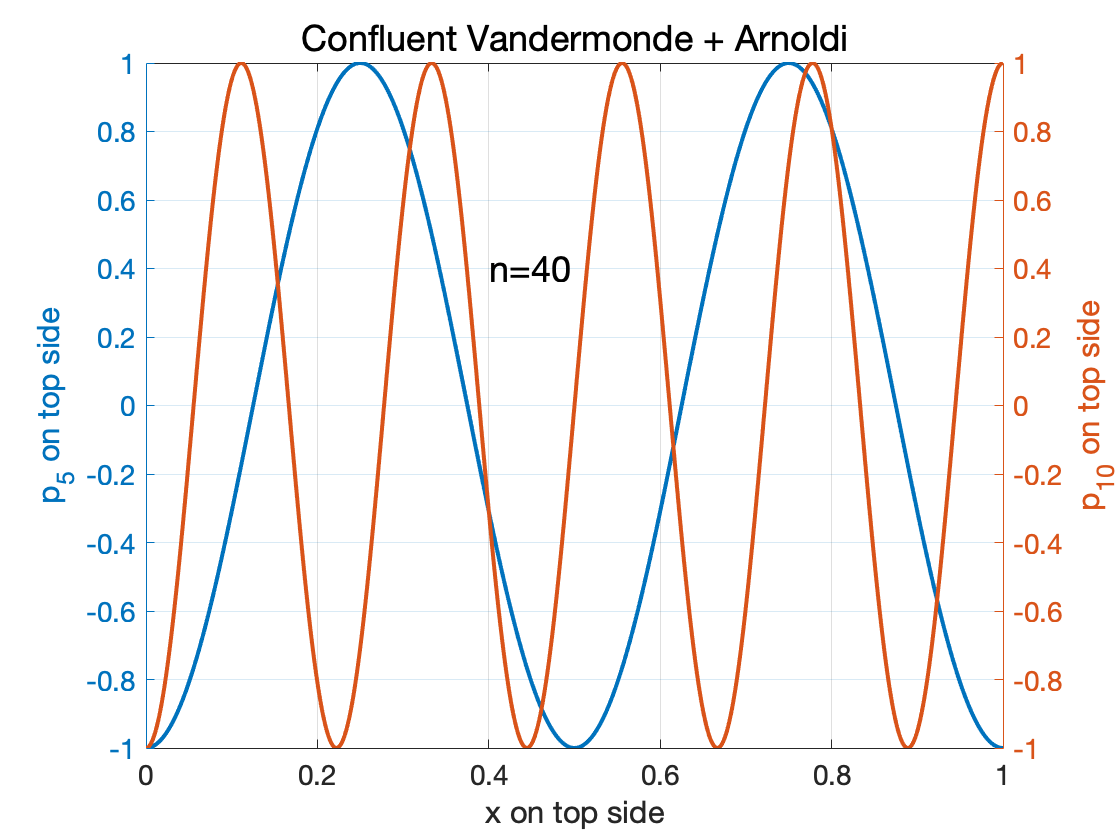}\\
  \includegraphics[scale=.16]{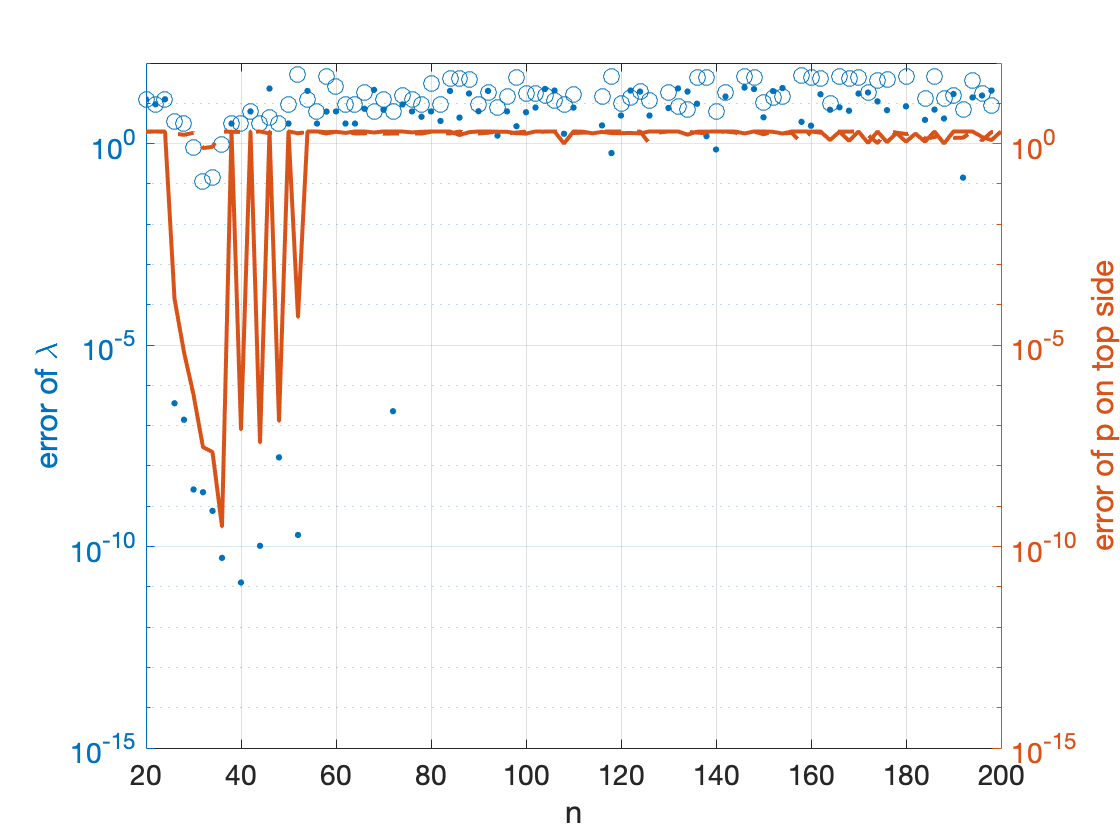}\quad%
  \includegraphics[scale=.16]{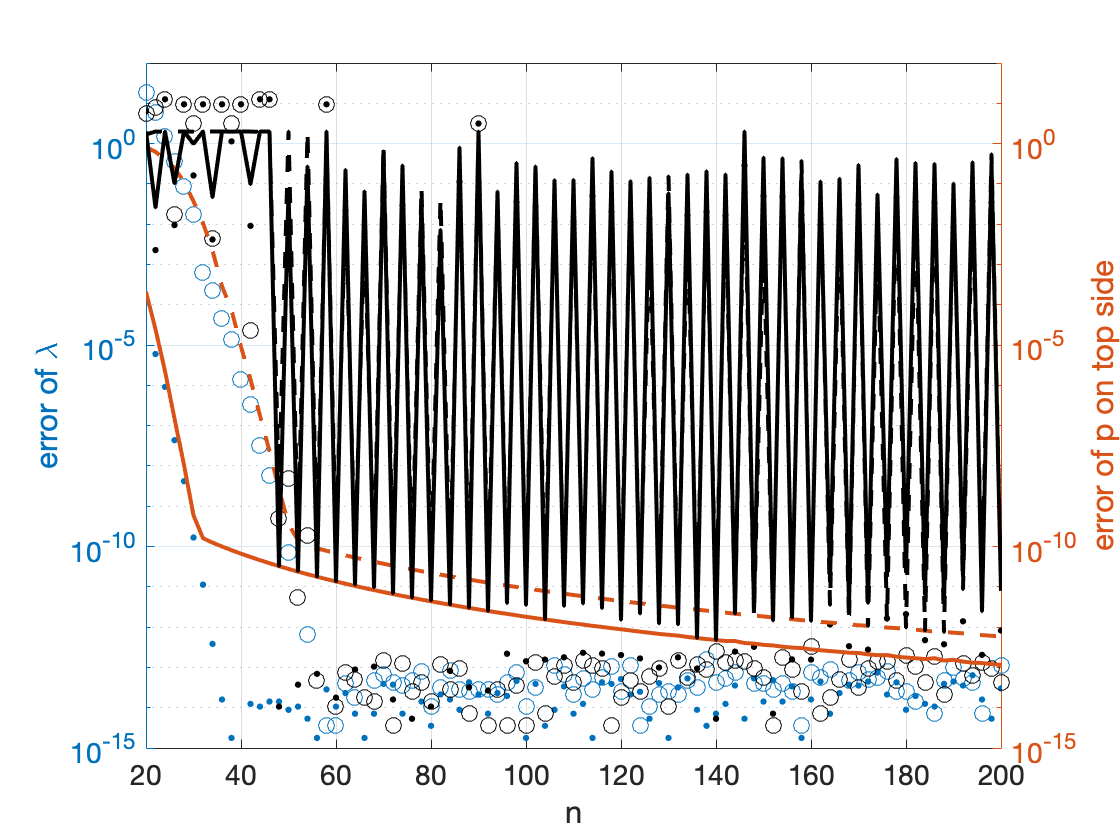}
  \caption{\small{Sloshing modes on the top side of a square cup by the confluent Vandermonde
      without (left column) and with (right column) Arnoldi. Top: the $k$th eigenfunction $p_k$
      ($\|p_k\|=1$), $k=5, 10$, for the eigenvalues $0=\lambda_1\le \ldots\le
      \lambda_k\le\ldots$. Bottom: errors of the 5th ($-$, $\cdot$) and 10th ($--$, $\circ$)
      eigenfunctions/eigenvalues. Bottom right: Hermite vs function values (in black) orthogonalized
      bases.}}
  \label{figex6}
\end{figure}

This example is taken from \cite{mora2015}. The domain is $\Omega=(0,1)^2$. On the left, right and
bottom sides, the homogeneous Neumann condition $\partial_{\boldsymbol\nu}u=0$ is imposed, while on
the top side is the eigenvalue equation $\partial_{\boldsymbol\nu}u=\lambda u$. The analytic
solution is $\lambda_n=n\pi\tanh(n\pi)$, $u_n=\cos(n\pi x)\sinh(n\pi y)$; see \cite{mora2015}. Let
$I: u\mapsto u$ be the identity operator on the top side and $O: u\mapsto 0u$ be the zero operator
on the other sides. Define $B$ on $\partial\Omega$ piecewise equal to $I$ and $O$. Then, we can
write the boundary equation uniformly as $\partial_{\boldsymbol\nu}u=\lambda Bu$ on
$\partial\Omega$. Using the notation from the previous section, the discrete rectangular eigenvalue
problem becomes $Q_{\nu}\boldsymbol{\beta}=\lambda BQ_{00}\boldsymbol{\beta}$.  Now, testing the
equation with the column space basis of $BQ_{00}$ as before is not a good idea because $BQ_{00}$ is
rank deficient. Note that the two sides of the equation live in the sum of the column spaces of
$Q_{\nu}$ and $BQ_{00}$. So, we follow \cite{hashemi2021} to do the economic SVD:
$(Q_{\nu}, BQ_{00}) = USV^*$, but we keep the first $2n+1$ columns of $U$ as $\hat{U}$ to match the
length of $\boldsymbol{\beta}$, then multiply the eigenvalue equation with $\hat{U}^*$ to get
$\hat{U}^*Q_{\nu}\boldsymbol{\beta}= \hat{U}^*BQ_{00}\boldsymbol{\beta}$. The results are shown in
Figure~\ref{figex6}. On each side of the square, we used $20(n+1)$ Chebyshev points of the first
kind for discretization of the boundary equation. From the bottom right, we can see that the
function values orthogonalized basis leads to instability. Convergence of the eigenfunctions from
the method with Hermite Arnoldi basis slowers down after some $n$, which could be fixed with more
sampling points. For example, for $n=60$ and $100(n+1)$ points per side, we got the errors
$2.4\times10^{-14}$, $1.2\times 10^{-13}$ for $p_5$, $p_{10}$.

\section{Example 7: indefinite integral on two intervals.}

\begin{figure}[h]
  \centering
  \includegraphics[scale=.16]{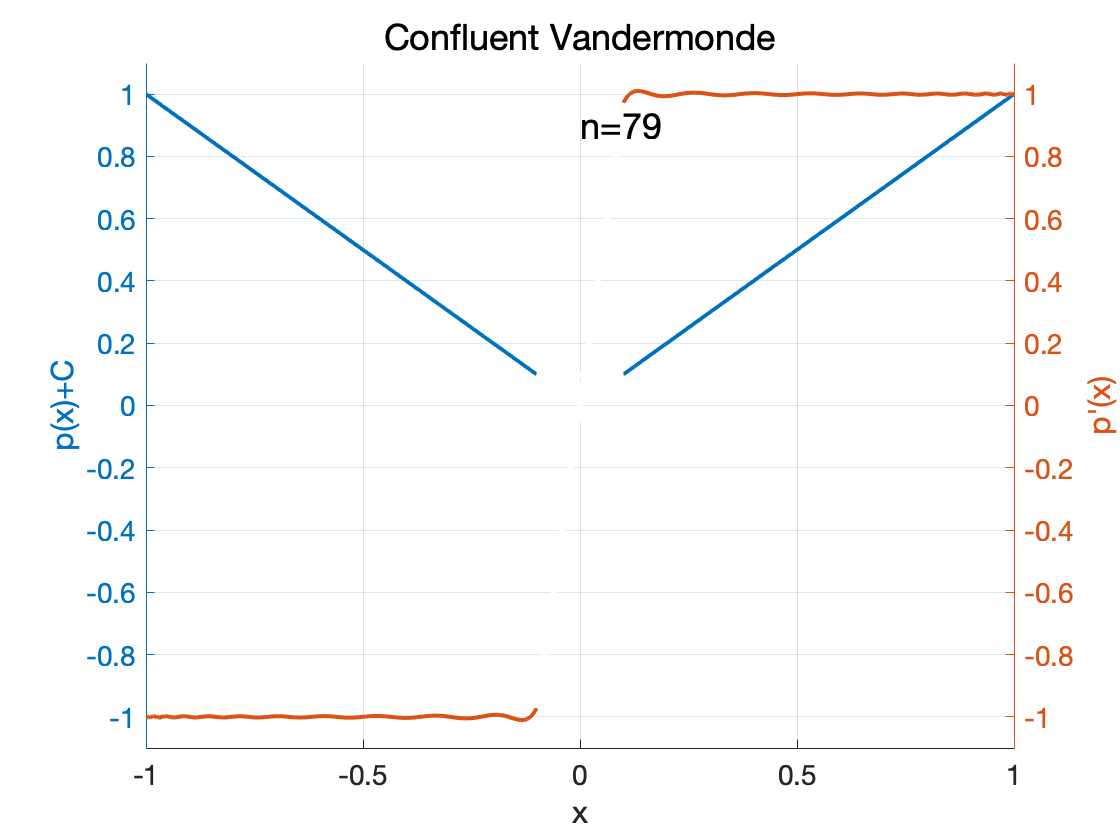}\quad%
  \includegraphics[scale=.16]{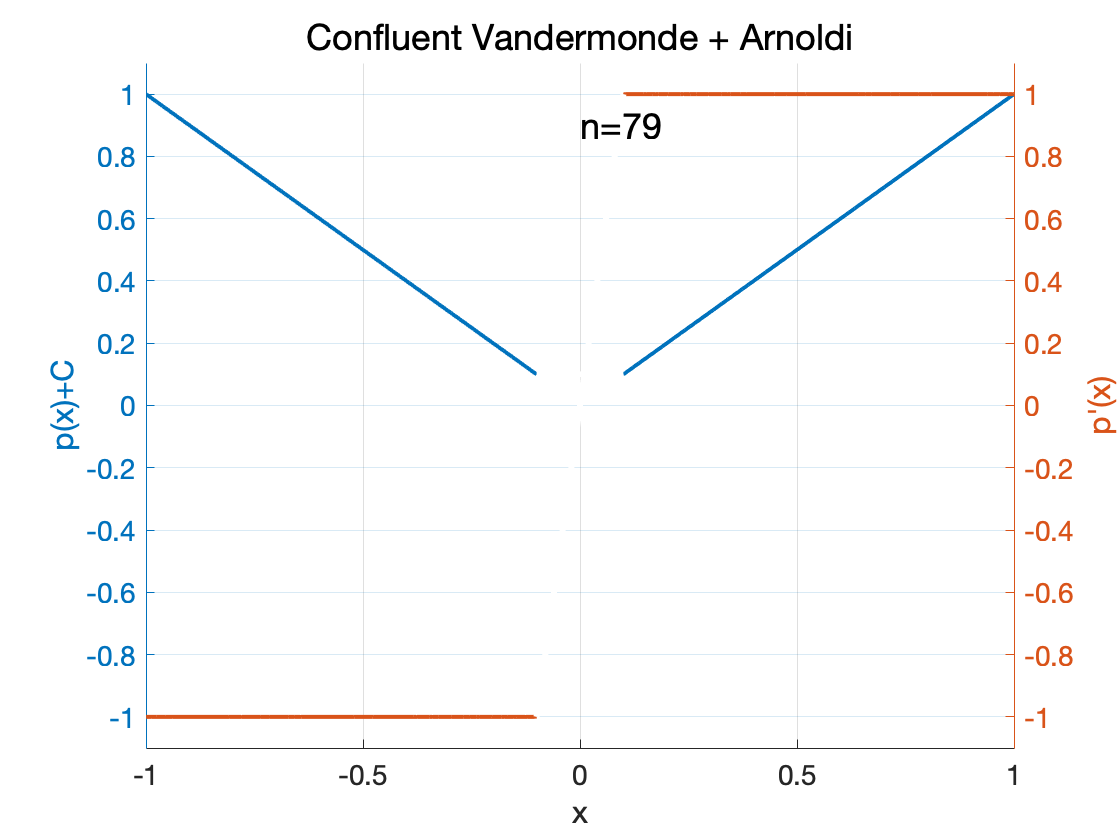}\\
  \includegraphics[scale=.16]{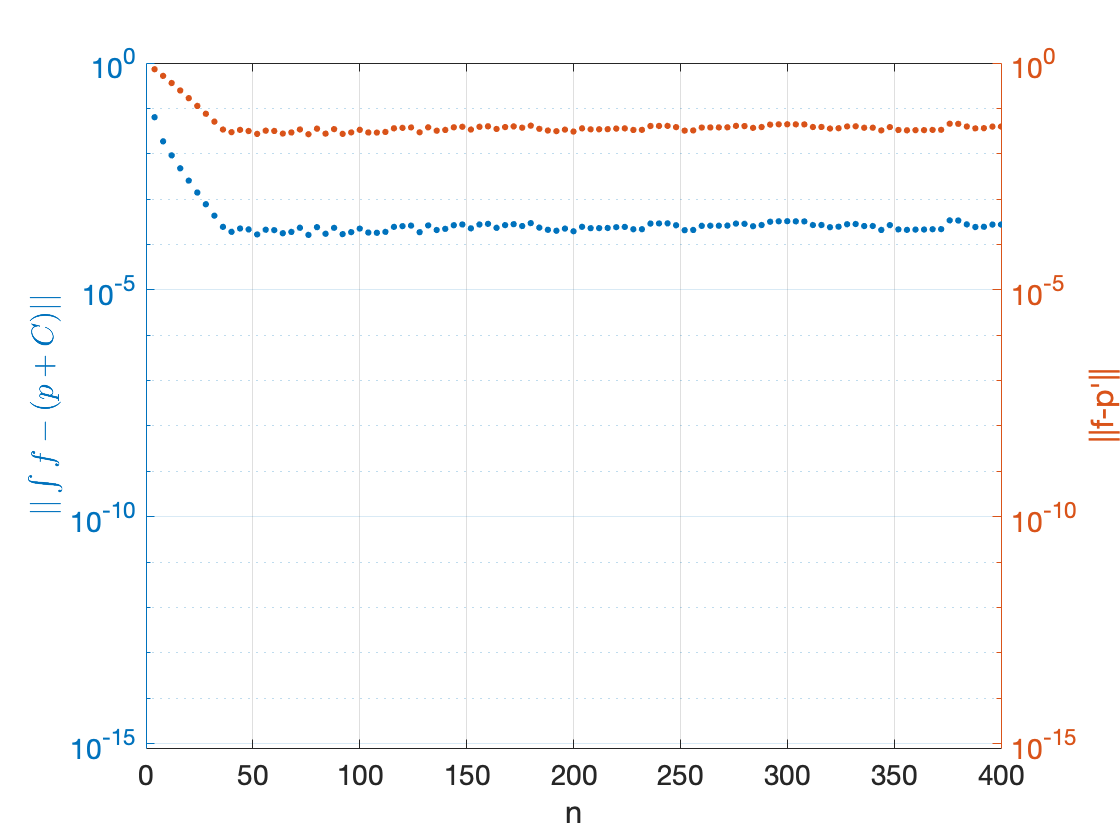}\quad%
  \includegraphics[scale=.16]{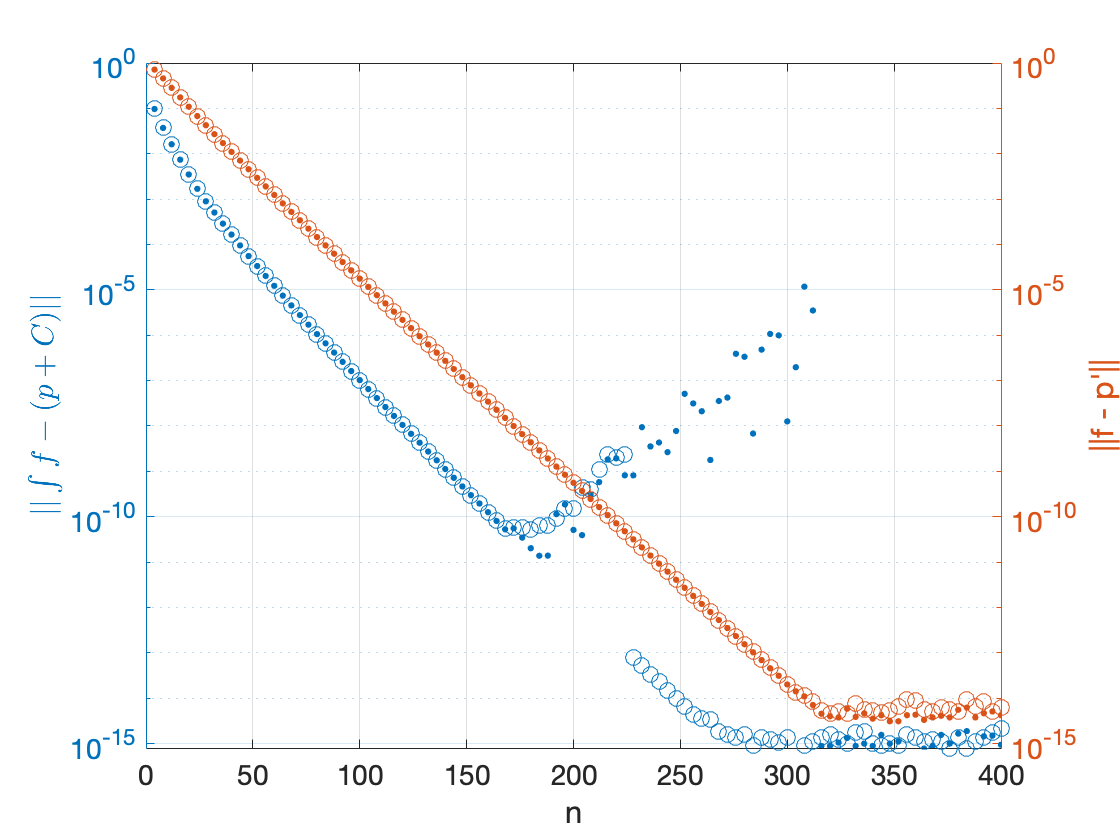}
  \caption{\small{Polynomial indefinite integral of $f(x)=\mathrm{sign}\ x$ on equidistant grids of
      $[-1,-1/10]\cup[1/10,1]$ by the confluent Vandermonde without (left column) and with (right
      column) Arnoldi. Bottom right: Hermite (dots) vs function values (circles) orthogonalized
      bases.}}
  \label{figex7}
\end{figure}

We want to find a polynomial approximation of $\int f(x)\mathrm{d}x\approx p(x)+C$, $C$ is a general
constant. For example, if $g(x)=\int_a^xu(t)~\mathrm{d}t$, we have $g'(x)= f(x)$, $g(a)=0$, the
simplest differential equation with a point value condition. The goal is to find $p(x)\approx g(x)$.
The implementation is mostly the same as \verb|polyfitAh| but the function header and the last line:
\verb|function [d, H] = polyfitAih(x,f,n)| and \verb|d = Q(m+1:2*m,2:n+1)\f; d = [0; d];|
alternatively, one can orthogonalize only the function values and from the resulting Hessenberg
matrix generate the derivative matrix. The solved example is shown in Figure~\ref{figex7}.  From the
bottom right, we can see that, for both the Hermite Arnoldi basis and the function values Arnoldi
basis, the error of $p$ can grow with $n$ in some range but still decreases to $10^{-15}$
afterwards. Such transient instability does not appear for $f(x)=1/(1+25x^2)$ on the Chebyshev grids
of $[-1,1]$.

\section{Discussion.}

This work is greatly inspired by the fascinating paper {\it Vandermonde with Arnoldi} \cite{va2021}
by Brubeck, Nakatsukasa and Trefethen. We saw a door open to the well-conditioned bases of
polynomial spaces on arbitrary grids, and to many interesting applications, e.g., the lightning
solvers \cite{lt2022, aa2021, log2021, trefethen2020}. In the spirit of Chebfun \cite{chebfun}, it
is natural then to consider differential, integral and other operations of polynomials under the
Arnoldi basis. Under the standard basis, the differential operation or Hermite interpolation leads
to the confluent Vandermonde matrix, which is ill-conditioned \cite{higham1990}. The Arnoldi process
can come to the rescue. However, it is not so straightforward to apply the method from
\cite{va2021}. The key idea behind the new method is the discovery of a connection of the confluent
Vandermonde matrix to a Krylov subspace. This enables us to represent and compute the derivatives in
a stable way. Two approaches are proposed to accomplish this task. One is constructing the basis for
the column space of the confluent Vandermonde matrix directly by the {\it Confluent Vandermonde with
  Arnoldi} method based on the confluent Krylov subspace that we discovered.  Another approach is
using the basis generated by {\it Vandermonde with Arnoldi} and constructing the derivative matrix
of the basis by exploiting the confluent Krylov subspace.  With the new functionality of computing
derivatives, we are able to do the Hermite interpolation/evaluation to high precision, to accomodate
the Neumann/Robin boundary conditions, to compute the Dirichlet-to-Neumann map and its eigenmodes,
and to do the indefinite integral.

We realized just upon finalizing this note that a formula of the derivative matrix has been given in
\cite{log2021, lt2022}. It is based on the recursive relation of the polynomial basis from {\it
  Vandermonde with Arnoldi}. Recall that the Arnoldi process produces the basis matrix $Q$ and the
Hessenberg matrix $H$ that satisfy $XQ_{-}=QH$. If we imagine that the number of rows of $X$ and $Q$
goes to infinity, then we would get the polynomial recursion
$ xp_k(x) = \sum_{j=1}^kp_j(x)h_{jk} + p_{k+1}(x)h_{k+1,k} $ or
$ p_{k+1}(x) = (xp_k(x) - \sum_{j=1}^kp_j(x)h_{jk}) / h_{k+1,k}$.  Taking derivative gives
$p_{k+1}'(x) = (p_k(x) + xp_k'(x) - \sum_{j=1}^kp_j'(x)h_{jk}) / h_{k+1,k}$. It can be seen that our
second approach (e.g. \verb|polyfitA| + \verb|polyvalAh|) is equivalent to this formula.

Comparison of the Hermite orthogonalized basis and the function values orthogonalized basis has been
made on each example. It is found that they perform equally well for the Hermite
interpolation/evaluation task. But the function values orthogonalized basis is unstable for the
example of sloshing modes, while the Hermite orthogonalized basis works stably. We also observed
instability of both the approaches for the Fourier extension example and the indefinite integral
example. A theoretical understanding of the comparison results and the instability phenomena could
be interesting.

\bibliographystyle{plain}
\bibliography{vandref.bib}

\end{document}